\documentclass[12pt,halfline,a4paper]{ouparticle}

\usepackage{graphicx}
\usepackage{natbib}
\usepackage[dvipsnames]{xcolor}
\usepackage{amsmath, amsthm, amssymb}
\usepackage{mathtools}
\usepackage{pdfpages}
\usepackage{afterpage}
\usepackage{hyperref}

\makeatletter
\renewcommand{\fnum@figure}{Fig. \thefigure}
\makeatother

\begin{document}
\title{Synthetic construction of the Hopf fibration\\ in a double orthogonal projection of \mbox{4-space}\\{\footnotesize This version is published in {\sl Journal of Computational Design and Engineering}, Volume {\bf 8}, Issue 3, June~2021, Pages 836–854, \texttt{https://doi.org/10.1093/jcde/qwab018}}}

\author{%
\name{Michal~Zamboj}
\address{Charles University, Faculty of Education, Department of Mathematics and Mathematical Education\\
M. D. Rettigové 4, 116 39 Prague 1, Czech Republic}
\email{michal.zamboj@pedf.cuni.cz}\\
}

\abstract{
The Hopf fibration mapping circles on a \mbox{3-sphere} to points on a \mbox{2-sphere} is well known to topologists. While the \mbox{2-sphere} is embedded in \mbox{3-space}, 4-dimensional Euclidean space is needed to visualize the \mbox{3-sphere}. Visualizing objects in 4-space using computer graphics based on their analytic representations has become popular in recent decades.
For purely synthetic constructions, we apply the recently introduced method of visualization of \mbox{4-space} by its double orthogonal projection onto
two mutually perpendicular \mbox{3-spaces} to investigate the Hopf fibration as a four-dimensional relation without  analogy in lower dimensions. In this paper, the method of double orthogonal projection is used for a direct synthetic construction of the fibers of a \mbox{3-sphere} from the corresponding points on a \mbox{2-sphere}. The fibers of great circles on the \mbox{2-sphere} create nested tori visualized in a stereographic projection onto the modeling \mbox{3-space}. The step-by-step construction is supplemented by dynamic 3D models showing simultaneously the \mbox{3-sphere}, \mbox{2-sphere}, and stereographic images of the fibers and mutual interrelations. Each step of the synthetic construction is supported by its analytic representation to highlight connections between the two interpretations.
}
\keywords{Hopf fibration; Hopf tori; four-dimensional visualization; stereographic projection; synthetic construction}
\date{\today}

\maketitle

\section{Introduction}
\label{sec:introduction}

\begin{figure}[!htb]
\centering
\includegraphics[width=0.6\linewidth]{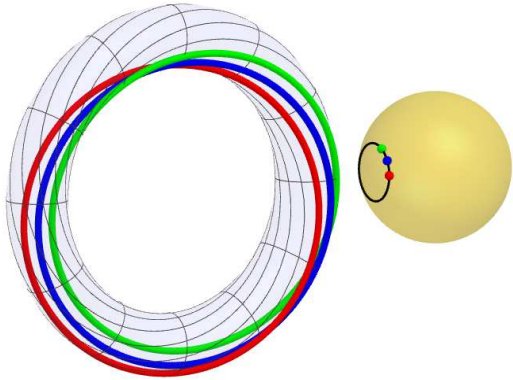}
\caption{The stereographic images of circular fibers corresponding to points on a circle on the base 2-sphere $\mathcal{B}^2$ are linked Villarceau circles on a torus.}
\label{fig:hopf_intro_torus}
\end{figure}

Mathematical visualization is an important instrument 
for understanding mathematical concepts. While analytic representations are convenient for proofs and analyses of properties, visualizations are essential for intuitive exploration and hypothesis making. Four-dimensional mathematical objects 
may lie beyond the reach of our three-dimensional imagination, but this is 
 not an obstacle to their mathematical description and study. Furthermore, using the modeling tools of computer graphics, we are able to construct image representations of four-dimensional objects to enhance their broader understanding. While many higher-dimensional 
mathematical objects are natural generalizations of lower-dimensional ones, the object of our study -- the Hopf fibration -- does not have this property. 
\begin{figure}[!htb]
\centering
\includegraphics[width=0.6\linewidth]{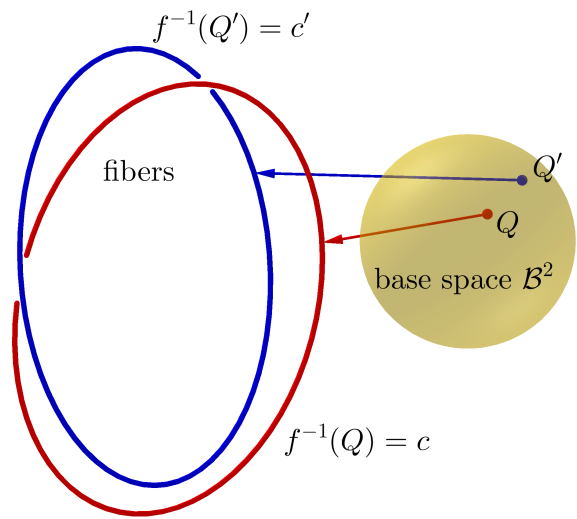}
\caption{Illustration of the Hopf fibration. Two distinct points $Q$ and $Q'$ on the 2-sphere $\mathcal{B}^2$ correspond to circular fibers $c$ and $c'$.}
\label{fig:hopf_intro}
\end{figure}
The Hopf fibration, introduced by \cite{Hopf1931a, Hopf1935}, defines a mapping between spheres of different dimensions. In this paper, we restrict ourselves to the correspondence it gives between spheres embedded in four- and three-dimensional spaces. To each point on a 2-sphere in three-dimensional space is assigned a circular fiber on a 3-sphere in four-dimensional space. The standard method of visualizing the fibers on the 3-sphere is to project them onto a three-dimensional space via stereographic projection, by which we can also grasp the topological properties of the Hopf fibration. Two distinct points on the 2-sphere correspond to disjoint circular fibers on the 3-sphere and their stereographic images are linked circles (Fig.~\ref{fig:hopf_intro}). To the points of a circle on the 2-sphere correspond circles on the 3-sphere that form a torus (Fig.~\ref{fig:hopf_intro_torus}).

Another facet of our interpretation 
 is the method of visualization itself. Instead of visualizing an analytic representation of a given object, we create the object directly in a graphical environment constructively point by point. To do so, we use generalized techniques of classical descriptive geometry to create three-dimensional images of four-dimensional objects. Our constructions or models being technically impossible to create by hand without computer graphics, we use interactive (or dynamic) 3D geometry software in which the images are placed in a virtual three-dimensional modeling space, and the observer is able to reach any point of this space and rotate views by manipulating the viewpoint in the graphical interface. Furthermore, in the models we shall present, the user can interactively manipulate (or animate) some parts or parameters of the model and immediately observe changes in the dependent objects.\footnote{Further on, by interactive models or visualizations we mean models with dynamically manipulable elements. This should not be confused with time-varying interactivity, often used in four-dimensional visualization.}

Two (not necessarily disjoint) groups might benefit from our paper: computer scientists interested in computer graphics, for whom we present a method of constructive visualization applied to the example of the Hopf fibration; and mathematicians, or students of mathematics, for whom we give a synthetic construction of the Hopf fibration and interactive graphical models to help better understand and intuitively explore it.

\subsection{Related work} 
Construction of the Hopf fibration is usually performed in abstract algebraic language as it involves the \mbox{3-sphere} embedded in $\mathbb{R}^4$ or $\mathbb{C}^2$.  However, recent research in computational geometry and graphics based on analytic representations has made possible partial video animations of the Hopf fibration along with models of stereographic images and visualizations in various software. The front covers of {\sl Mathematical Intelligencer}, vol. 8, no. 3 and vol. 9, no. 1 featured stereographic images of the Hopf fibration by \cite{Kocak1987}, taken from their pioneering computer-generated film projects with Banchoff, Bisshopp, and Margolis. \cite{Banchoff1990} wrote an comprehensive illustrated book on the fourth dimension, and its front cover has another inspiring picture of the Hopf fibration from his film. Moreover, moving toward 
the computer visualization of the Hopf fibration, \cite{Banchoff1988} constructed stereographic images of Pinkall's tori of given conformal type. The primary inspiration for writing the present paper was another film -- `Dimensions' by \cite{Dimensions}, in which the Hopf fibration is well explained and visualized in a variety of separate 
models. Coincidentally, another front cover -- of {\sl Notices of the AMS}, vol. 44, no. 5 -- was inspired by explanatory illustrations created in {\sl Wolfram Mathematica} published by \cite{Kreminski1997} in the context of the structure of the projective extension of real 3-space, $\mathbb{RP}^3$. A popular visualization of the Hopf fibration showing points on the base 2-sphere and the corresponding stereographic images of the fibers was created by \cite{Johnson2011} in the mathematics software {\sl Sage}. Johnson's code was modified by \cite{Chinyere2012} to visualize a similar fibration with trefoil knots as fibers instead of circles.

The visual aspect is foregrounded by \cite{Hanson2006a} (pp. 80-85 and 386-392) in visualizing quaternions on a \mbox{3-sphere}, and the author also describes them using Hopf fibrations. Visualizations of Hopf fibrations have been beneficial in topology in relation to the Heegaard splitting of a \mbox{3-sphere} by \cite{Canlubo2017}, in the use of quaternions in physics by \cite{OSullivan2015}, and in the description of motion in robotics by \cite{Yershova2010}. Based on analytic representations, \cite{Black2010} in chapter 6 of his dissertation depicts similar orthogonal projections of tori as do we, and supplemented them with animations.

Interactive tools for the visualization of two-dimensional images of four-dimensional objects from different viewpoints in 4-space were developed and applied by \cite{Zhou1991}. In another early-stage thesis in the area, \cite{Heng1992} (supervised by Hanson) wrote on the use of interactive mathematical visualization techniques in computer graphics applied to the exploration of 3-manifolds. The continuous development of interactive frameworks and methods of four-dimensional visualization is also apparent in subsequent work co-authored by Hanson (e.g. \cite{Hanson1999, Zhang2007,Thakur2007, Chu2009}, the last two papers including visualizations of a flat torus embedded in 4D from different viewpoints).

In this exposition, we use the language and elementary constructions of the double orthogonal projection described in \cite{Zamboj2018a}, and the constructions of sections of four-dimensional polytopes, cones, and spheres published in the series of articles \cite{Zamboj2018b,Zamboj2019a,Zamboj2019b}.

\subsection{Contribution} 
We present an application of a method of visualization using computer graphics for constructing and examining a phenomenon of 4-space that has no analogy in lower dimensions. We use the method of double orthogonal projection to create a purely synthetic graphical construction of the Hopf fibration. This paper contributes to the field in two directions: a novel graphical construction of the Hopf fibration and an application of the double orthogonal projection. In contrast to previous work on visualization of the Hopf fibration, in which separate illustrations created from analytic representations were graphical results or explanatory additions, we use mathematical visualization via the double orthogonal projection as a tool to synthetically construct the fibration. For this purpose, we revisit the 
analytic definition of the Hopf fibration (points on a circle on a \mbox{3-sphere} map to a point on a \mbox{2-sphere}, given by expression~(\ref{eq:hopfmap})), which gives us the solved puzzle, and initiate our synthetic approach by decomposing the fibration into pieces. After this, we propose an elementary step-by-step construction of the inverse process. From points on the \mbox{2-sphere}, we construct fibers of the \mbox{3-sphere} with the use of only elementary (constructive) geometric tools. On top of that, we construct the resulting stereographic images of the fibers in one complex graphical interpretation. Even though stereographic images are common in visualization of the Hopf fibration, the difference here lies in our synthetic construction of them in order to confirm our results and observations. Finally, the method of visualization is applied to provide a graphical analysis of the properties of cyclic surfaces on a 3-sphere, four-dimensional modulations, and filament packings. Our constructions (Figs.~\ref{fig:pointionsphere}--\ref{fig:hopf2fibers}, \ref{fig:hopf_SBS}--\ref{fig:hopf_2tori}) are supplemented by interactive 3D models in {\sl GeoGebra 5} (see the online GeoGebra Book \cite{ZambojGGBBook2019}).\footnote{The software {\sl GeoGebra 5} is used due to its general accessibility, but any other 3D interactive (dynamic) geometry software may be used with the same outcomes. For overall understanding, we strongly recommend following the 3D models simultaneously with the text.} The final visualizations (Figs. \ref{fig:nested_tori_S}--\ref{fig:packingsful}) and videos (Suppl. Files 8 and 9) are created in {\sl Wolfram Mathematica 11} for better graphical results. Throughout the paper, we give the relevant analytic background to our synthetic visual approach to enrich overall understanding of the Hopf fibration.

\subsection{Paper Organization}
In Section~\ref{sec:math} we introduce the analytic definition and basic properties of the Hopf fibration, which is followed in Section~\ref{sec:preliminary} by a brief description of how to depict a \mbox{3-sphere} in its double orthogonal projection onto two mutually perpendicular 3-spaces and how to construct its stereographic image. In the main part of the paper, Section~\ref{sec:construction}, we give a synthetic construction of a Hopf fiber on a \mbox{3-sphere}, corresponding to a point on a \mbox{2-sphere}, using only elementary tools. The resulting double orthogonal projection and stereographic images of tori on the \mbox{3-sphere}, corresponding to two families of circles on the \mbox{2-sphere}, are given in Section~\ref{sec:tori}. In Section~\ref{sec:applications} we apply the method to constructions of cyclic surfaces, visualization of four-dimensional modulations with respect to polyhedral arrangements of vertices on the 2-sphere, and filament packings of the 3-sphere. After concluding with perspectives on future work, an appendix gives the parametrizations used in the figures.

\section{Mathematical background}
\label{sec:math}

In algebraic topology, 
a fibration is a certain type of projection from one topological space (the total space) onto another (the base space) that decomposes the total space into fibers. We proceed to formally define these and other terms that are needed in the sequel.\\\\
A {\sl topological space} is an ordered pair $(X, \tau)$, where $X$ is a set and $\tau$ is a collection of subsets of $X$ satisfying the following conditions:
\begin{enumerate}
\item The empty set and $X$ belong to $\tau$.
\item Arbitrary unions of sets in $\tau$ belong to $\tau$.
\item The intersection of any finite number of sets in $\tau$ belongs to $\tau$.
\end{enumerate}
$\tau$ is called the {\sl topology} on $X$ of the topological space $(X,\tau)$.\\\\
Let $\mathcal{T}$ and $\mathcal{B}$ be topological spaces and $f: \mathcal{T} \rightarrow \mathcal{B}$ a continuous surjective mapping. The ordered triple $(\mathcal{T}, f ,\mathcal{B})$ is called the {\sl fiber space} or {\sl fibration}, $\mathcal{T}$ is the {\sl total space}, $\mathcal{B}$ is the {\sl base space}, and $f$ is the {\sl projection} (or {\sl fibration}\footnote{The term fibration is sometimes used for the whole triple, sometimes for the component mapping of the triple. For the Hopf fibration we freely use either convention throughout the paper as it should not cause any confusion.}) {\sl of the fibre space}. The inverse image of a point (element) of the base space $\mathcal{B}$ in the projection $f$ is the {\sl fiber above this point}. \\\\
A {\sl \mbox{3-sphere}} is the three-dimensional boundary of a \mbox{4-ball} in four-dimensional Euclidean space. It is a natural generalization of a \mbox{2-sphere} defined as the two-dimensional boundary of a ball in three-dimensional Euclidean space.\\

\begin{figure}[!htb]
\centering
\includegraphics[width=0.6\linewidth]{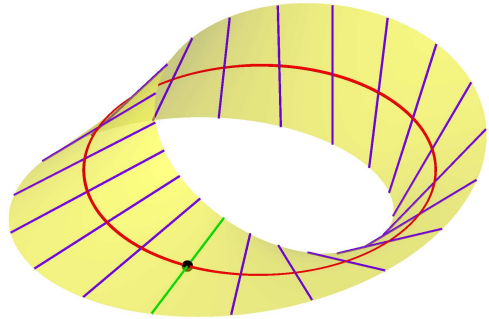}
\caption{Möbius band with central circle as the base space and segments as fibers. Each segment is projected onto its point of intersection with the central circle.}
\label{fig:mobius}
\end{figure}

Some elementary examples may help illustrate the above definition of a fibration. \cite{Hanson2006a}, p. 386, gives a shag rug as an informal example of a trivial fibration, where the backing of the rug is the base space and the threads are fibers.
Similarly, a spiky massage ball has the 2-sphere as the base space and its spikes are fibers above points of the sphere. The Hopf fibration, defined shortly, will also have the 2-sphere as the base space, the fibers will be circles (not necessarily glued to the base space), but the topology of the total space will be less trivial. An important nontrivial example of a fibration is a Möbius band as the total space (see also \cite{Shoemake1994}), with the central circle as the base space, and fibers being a collection of line segments twisting around the circle (see Fig.~\ref{fig:mobius}). The projection maps each point of the Möbius band, which belongs to exactly one of the line segments, to its point on the central circle, and the inverse image of a point on the central circle is the whole line segment. 
Finally and most importantly for us, the Hopf fibration is a mapping from a \mbox{3-sphere} (the total space $\mathcal{T}^3$, the exponent indicates the dimension of the space) to a \mbox{2-sphere} (the base space $\mathcal{B}^2$) such that each distinct circle (fiber) on the \mbox{3-sphere} $\mathcal{T}^3$ corresponds to a distinct point on the \mbox{2-sphere} $\mathcal{B}^2$. In the other direction, the total space $\mathcal{T}^3$ is formed by circular fibers corresponding to points in the base space (the \mbox{2-sphere} $\mathcal{B}^2$). Furthermore, inverse images of points on the \mbox{2-sphere} correspond to non-intersecting linked Villarceau circles on a torus in the \mbox{3-sphere} when stereographically projected (Fig.~\ref{fig:hopf_intro_torus}).

An analytic definition of the Hopf fibration in coordinate geometry (expression~(\ref{eq:hopfmap})) is given after we have made some preliminary remarks. For an elementary introduction to the Hopf fibration with visualizations, see \cite{Lyons2003, Treisman2009, Ozols2007}, and for further details and proofs in modern topological language, see the textbook by \cite{Hatcher2002}, Chapter~4.\\\\

Let  $\{X_i\}_{i\in I}$ be a collection of sets indexed by a set $I$. The set of functions $f: I\rightarrow \bigcup_{i\in I}X_i$ such that $f(i)\in X_i$ for each $i\in I$ is called the {\sl Cartesian} (or {\sl direct}) {\sl product} of the family of sets $\{X_i\}_{i\in I}$. The Cartesian product is denoted by $\prod_{i\in I}X_i$ or $X_1\times X_2\times\dots$.\\\\

For the purposes of visualization, we follow the direct construction of the Hopf fibration in real four-dimensional space, $\mathbb{R}^4$, given for example in \cite{Treisman2009}; we also give an alternative construction in complex space useful for simple calculations. As a valuable by-product we thus obtain visualizations of objects embedded not only in $\mathbb{R}^4$ but also in $\mathbb{C}^2$, where $\mathbb{C}^2=\mathbb{C}\times\mathbb{C}$ is the Cartesian product of two complex coordinate systems. 

\paragraph{Analytic definition of the Hopf fibration.} At first, we construct a mapping that projects a point on the \mbox{3-sphere} $\mathcal{T}^3$ to a point on the  \mbox{2-sphere} $\mathcal{B}^2$. The unit \mbox{2-sphere} $\mathcal{B}^2$ 
in $\mathbb{R}^3$ is the set of points given by the equation
\begin{equation}
\label{eq:2sphere}
x^2+y^2+z^2=1.
\end{equation}
The unit \mbox{3-sphere} $\mathcal{T}^3$ 
in $\mathbb{R}^4$ is given by the analogous equation
\begin{equation}
\label{eq:3sphere}
x^2+y^2+z^2+w^2=1.
\end{equation}
Given complex numbers $z_1=x_0+iy_0$ and $z_2=z_0+iw_0$, under the identification of $\mathbb{C}$ with $\mathbb{R}^2$ we have $z_1=[x_0,y_0]$ and $z_2=[z_0,w_0]$.\footnote{The labeling $x_0,y_0,z_0,w_0$ refers to the location of a point with respect to the coordinate axes $x,y,z,w$, which will be graphically interpreted later.} This way, the unit \mbox{3-sphere} $\mathcal{T}^3$ embedded in $\mathbb{C}^2$ has the equation
\begin{equation}
\label{eq:3sphereC}
|z_1|^2+|z_2|^2=1.
\end{equation}
Likewise, setting $\zeta=x+iy$, 
we can reformulate equation~(\ref{eq:2sphere}) as an embedding of $\mathcal{B}^2$ in $\mathbb{C}\times\mathbb{R}$:
\begin{equation}
\label{eq:2sphereC}
|\zeta|^2+z^2=1.
\end{equation}
The mapping
\begin{equation}
\label{eq:hopfmap}
\begin{split}
\displaystyle f:\mathbb{C}^2\to\mathbb{C}\times\mathbb{R}\text{ such that }\\
f(z_1,z_2)=(2z_1\overline{z_2},|z_1|^2-|z_2|^2)
\end{split}
\end{equation}
is called the {\sl Hopf fibration}.
\footnote{At this point, the Hopf fibration is given by its existential definition. Other equivalent definitions of the Hopf fibration as a ratio of complex numbers (see \cite{Hatcher2002}) or via quaternions (see \cite{Hanson2006a}) would require a much broader theoretical discourse and would not fit our perspective any better.}


\paragraph{Proposition 1.} 
The preimage $[z_1,z_2]$ of the point $[\zeta,z]$ on $\mathcal{B}^2\subset\mathbb{C}\times\mathbb{R}$ under the Hopf fibration (\ref{eq:hopfmap}) is on $\mathcal{T}^3\subset\mathbb{C}^2$.\\\\

Substituting $\zeta=2z_1\overline{z_2}$ and $z=|z_1|^2-|z_2|^2$ into equation~(\ref{eq:2sphereC}) yields
\begin{equation}
\label{eq:pointon3sphere}
\begin{split}
|2z_1\overline{z_2}|^2+(|z_1|^2-|z_2|^2)^2&=1\\
4((z_1\overline{z_2})(\overline{z_1}z_2))+|z_1|^4-2(z_1\overline{z_1})(z_2\overline{z_2})+|z_2|^4&=1\\(|z_1|^2+|z_2|^2)^2&=1\\\text{and since } |z_1|^2+|z_2|^2\geq 0 \text{, we have }|z_1|^2+|z_2|^2&=1.
\end{split}
\end{equation}
According to equation~(\ref{eq:3sphereC}), the point $[z_1,z_2]$ lies on $\mathcal{T}^3$.

Next, we describe how to create a circle from a point on $\mathcal{T}^3$ such that its image is a fixed point on $\mathcal{B}^2$.

\paragraph{Proposition 2.} The Hopf fibration (\ref{eq:hopfmap}) maps all points of a circular fiber on $\mathcal{T}^3\subset\mathbb{C}^2$ to one point on $\mathcal{B}^2\subset\mathbb{C}\times\mathbb{R}$. \\\\
Let $[z_1,z_2]$ be a point on the 3-sphere $\mathcal{T}^3\subset\mathbb{C}^2$. Let $\lambda=l_1+il_2\in\mathbb{C}$ be such that $|\lambda|^2=1$, i.e. $\lambda$ represents a point on a unit circle embedded in $\mathbb{C}$.  From equation~(\ref{eq:pointon3sphere}), for the point $[\lambda z_1,\lambda z_2]\in\mathbb{C}^2$ it holds that
\begin{equation}
|\lambda z_1|^2+|\lambda z_2|^2=|\lambda|^2(|z_1|^2+|z_2|^2)=|z_1|^2+|z_2|^2=1, \text{for each }\lambda.
\label{eq:prop2_1}
\end{equation}
Hence the point $[\lambda z_1,\lambda z_2]$ lies on $\mathcal{T}^3\subset\mathbb{C}^2$.\\

Rewriting the point $[\lambda z_1,\lambda z_2]\in\mathbb{C}^2$ in its parametric representation in $\mathbb{R}^4$ with $\lambda=l_1+il_2, z_1=x_0+iy_0$ and $z_2=z_0+iw_0$, we have
\begin{equation}
\label{eq:imre}
\begin{split}
\begin{pmatrix}
	\operatorname{Re}(\lambda z_1)\\
	\operatorname{Im}(\lambda z_1)\\
	\operatorname{Re}(\lambda z_2)\\
	\operatorname{Im}(\lambda z_2)\\
\end{pmatrix}=
\begin{pmatrix}
	l_1x_0-l_2y_0\\
	l_1y_0+l_2x_0\\	
	l_1z_0-l_2w_0\\
	l_1w_0+l_2z_0\\
\end{pmatrix}.
\end{split}
\end{equation}
For each $[l_1,l_2]\in\mathbb{R}^2$, the last expression defines a set of points $l_1\overrightarrow{u}+l_2\overrightarrow{v}$ in a plane in $\mathbb{R}^4$ through the origin $O=[0,0,0,0]$ with the directional vectors $\overrightarrow{u}=(x_0,y_0,z_0,w_0)$ and $\overrightarrow{v}=(-y_0,x_0,-w_0,z_0)$. By equation~(\ref{eq:prop2_1}), the set of all points $[\lambda z_1,\lambda z_2]$ for each $\lambda$ is the intersection of $\mathcal{T}^3$ and a plane through its center. Hence the set of all points $[\lambda z_1,\lambda z_2]$ for $\lambda\in\mathbb{C}$ is a unit circle on $\mathcal{T}^3$.\\
By equation~(\ref{eq:hopfmap}) defining the Hopf fibration $f$, 
\begin{equation}
\begin{split}
f(\lambda z_1,\lambda z_2)=(2\lambda z_1\overline{\lambda z_2},|\lambda z_1|^2-|\lambda z_2|^2)=(2|\lambda|^2 z_1\overline{z_2},|\lambda|^2(|z_1|^2-|z_2|^2))\\=(2z_1\overline{z_2},|z_1|^2-|z_2|^2)=f(z_1, z_2).
\end{split}
\end{equation}
so that $f$ maps $[\lambda z_1, \lambda z_2]$ for $|\lambda|^2=1$ to the same point on the unit \mbox{2-sphere} $\mathcal{B}^2$ as $[z_1,z_2]$.

{\paragraph{Hopf coordinates.} In our synthetic reconstruction, we construct the inverse mapping -- for a point on $\mathcal{B}^2$, we find the circle on $\mathcal{T}^3$. The point on $\mathcal{B}^2$ will be constructed in terms of its polar and azimuthal angles in the representation of $\mathcal{B}^2$ in spherical coordinates. We thus seek a relation between the spherical coordinates of a point on $\mathcal{B}^2$ and a point $[z_1,z_2]$ on $\mathcal{T}^3$ in its trigonometric representation $z_1=r_A(\cos\alpha+i\sin\alpha)$ and $z_2=r_B(\cos\beta+i\sin\beta)$ for $r_A,r_B\geq0$ and $\alpha,\beta\in\mathbb{R}$. From equation~(\ref{eq:pointon3sphere}), for a point $[z_1,z_2]$ on $\mathcal{T}^3$ it holds that $|z_1|^2+|z_2|^2=r_A^2+r_B^2=1$, and so there exists a unique $\gamma\in\langle0,\frac{\pi}{2}\rangle$ such that $r_A=\cos\gamma$ and $r_B=\sin\gamma$. Then a point on $\mathcal{T}^3$ has the following coordinates in~$\mathbb{R}^4$:
\begin{equation}
\label{eq:parS3gon}
\begin{pmatrix}
\cos\gamma\cos\alpha\\
\cos\gamma\sin\alpha\\
\sin\gamma\cos\beta\\
\sin\gamma\sin\beta
\end{pmatrix}, \alpha, \beta \in \langle0,2\pi), \gamma \in \Big\langle0,\frac{\pi}{2}\Big\rangle.
\end{equation}
With the use of this representation, the image $f(z_1,z_2)$ of the point $[z_1,z_2]$ under the Hopf fibration (\ref{eq:hopfmap}) has first coordinate
\begin{equation}
\label{eq:1stcoordinate}
\begin{split}
2z_1\overline{z_2}=2\cos\gamma(\cos\alpha+i\sin\alpha)\sin\gamma(\cos\beta-i\sin\beta)=\\
2\cos\gamma\sin\gamma(\cos\alpha\cos\beta+\sin\alpha\sin\beta+\\i(\sin\alpha\cos\beta-\cos\alpha\sin\beta))=
\sin(2\gamma)(\cos(\alpha-\beta)+i\sin(\alpha-\beta)),
\end{split}
\end{equation}
and second coordinate
\begin{equation}
\label{eq:2ndcoordinate}
\begin{split}
|z_1|^2-|z_2|^2=r_A^2-r_B^2=\cos^2\gamma-\sin^2\gamma=\cos(2\gamma).
\end{split}
\end{equation}
Let $2\gamma=\psi\in\langle0,\pi\rangle$ and $\alpha-\beta=\varphi\in\langle0,2\pi)$ (taken modulo $2\pi$, and likewise for all further operations with angles). Considering the real and imaginary parts in (\ref{eq:1stcoordinate}
) and the third coordinate in (\ref{eq:2ndcoordinate}), the coordinates of the image $f(z_1,z_2)$ in $\mathbb{R}^3$ are
\begin{equation}
\begin{split}
\label{eq:spherical}
\begin{pmatrix}
\sin(2\gamma)\cos(\alpha-\beta)\\
\sin(2\gamma)\sin(\alpha-\beta)\\
\cos(2\gamma)
\end{pmatrix}=
\begin{pmatrix}
\sin\psi\cos\varphi\\
\sin\psi\sin\varphi\\
\cos\psi
\end{pmatrix},\\ \psi\in\langle0,\pi\rangle,\varphi\in\langle0,2\pi).
\end{split}
\end{equation}
We thus obtain the spherical coordinates of the image of a point on $\mathcal{T}^3$, which lies on the unit \mbox{2-sphere} $\mathcal{B}^2$, and we have also reduced the number of parameters (from three $\alpha,\beta,\gamma$ to two $\varphi,\psi$). Finally, to construct the preimage of a point on $\mathcal{B}^2$, 
we substitute $\varphi$ and $\psi$ into equation (\ref{eq:parS3gon}) parametrizing $\mathcal{T}^3$ 
in $\mathbb{R}^4$, thereby obtaining the {\sl Hopf coordinates of a \mbox{3-sphere}}:
\begin{equation}
\label{eq:hopfcoordinates}
\begin{pmatrix}
\cos\frac{\psi}{2}\cos(\varphi+\beta)\\
\cos\frac{\psi}{2}\sin(\varphi+\beta)\\
\sin\frac{\psi}{2}\cos\beta\\
\sin\frac{\psi}{2}\sin\beta
\end{pmatrix}, \beta,\varphi \in \langle0,2\pi),\psi\in\langle0,\pi\rangle.
\end{equation}
Represented in $\mathbb{C}^2$, we have
\begin{equation}
\label{eq:hopfcoordinatescomplex}
\begin{pmatrix}
z_1\\
z_2
\end{pmatrix}=
\begin{pmatrix}
\cos\frac{\psi}{2}(\cos(\varphi+\beta)+i\sin(\varphi+\beta))\\
\sin\frac{\psi}{2}(\cos\beta+i\sin\beta)
\end{pmatrix}.
\end{equation}

Let us show that two fibers are disjoint (see also \cite{Treisman2009}).}

\paragraph{Proposition 3.} Hopf fibers are disjoint circles.\\\\
Let $\mathcal{T}^3$ be a 3-sphere given by equation~(\ref{eq:3sphere}). Its intersection with the \mbox{3-space} $w=0$ is the (equatorial) \mbox{2-sphere} $\mathcal{B}^2$ with the equation~(\ref{eq:2sphere}).

Consider a point $A$ on $\mathcal{T}^3$ with coordinates $A[1,0,0,0]$ in $\mathbb{R}^4$. The set of points $c_A=\lambda_A A$ for $|\lambda_A|^2=1,\lambda_A\in\mathbb{C}$, is the unit circle defined by the rotation of $A$ by $\lambda_A$ about the origin in the plane $(x,y)$. Now let $B[x_B,y_B,z_B,w_B]$ be another point on $\mathcal{T}^3$ and not on the circle $c_A$ through $A$ so that $(z_B,w_B)\neq(0,0)$. The set of points $c_B=\lambda_B B$ for $|\lambda_B|^2=1,\lambda_B\in\mathbb{C}$ is by equation~(\ref{eq:prop2_1}) a unit circle on $\mathcal{T}^3$ with parametric representation 
\begin{equation}
\label{eq:Bimre}
\begin{split}
\begin{pmatrix}
	\operatorname{Re}(\lambda_B \sqrt{x_B^2+y_B^2} (x_B+iy_B))\\
	\operatorname{Im}(\lambda_B \sqrt{x_B^2+y_B^2}(x_B+iy_B))\\
	\operatorname{Re}(\lambda_B \sqrt{z_B^2+w_B^2}(z_B+iw_B))\\
	\operatorname{Im}(\lambda_B \sqrt{z_B^2+w_B^2}(z_B+iw_B))\\
\end{pmatrix}.
\end{split}
\end{equation}
The unit circle $c_B$ intersects the equatorial unit \mbox{2-sphere} $\mathcal{B}^2$ only if some point on $c_B$ is in the \mbox{3-space} $w=0$. Let us again use the trigonometric representation to show that there are only 2 intersections of $c_B$ with $\mathcal{B}^2$ (just as for great circles on a \mbox{2-sphere} intersecting its equator). Then $z_B + iw_B=\sqrt{z_B^2+w_B^2}(\cos\delta+i\sin\delta)$ for $\delta \in \langle0,2\pi)$, and $\lambda_B=\cos\varkappa+i\sin\varkappa$ for $\varkappa\in \langle0,2\pi)$. For the point on $c_B$ in the 3-space $w=0$, the equation 
\begin{equation}
\begin{split}
\operatorname{Im}(\lambda_B \sqrt{z_B^2+w_B^2}(z_B+iw_B)))=\operatorname{Im}((\cos\varkappa+i\sin\varkappa)(\sqrt{z_B^2+w_B^2}(\cos\delta+i\sin\delta)))\\=\sqrt{z_B^2+w_B^2}\sin(\varkappa+\delta)=0
\end{split}
\end{equation}
has, for $(z_B,w_B)\neq(0,0)$, two solutions for $\lambda_B$: $\varkappa=-\delta \mod{2\pi}$ or $\varkappa=\pi-\delta \mod{2\pi}$ for each $\delta \in \langle0,2\pi)$ . Therefore, $c_B$ has only two antipodal points $K_1$ and $K_2$ in common with $\mathcal{B}^2$. 

Since for $z$-coordinates of $c_B$ we have
\begin{equation}
\begin{split}
\operatorname{Re}((\cos\varkappa+i\sin\varkappa)(\sqrt{z_B^2+w_B^2}(\cos\delta+i\sin\delta)))=\sqrt{z_B^2+w_B^2}\cos(\varkappa+\delta)\neq0\\\text{ for }\varkappa=-\delta\text{ or }\varkappa=\pi-\delta,   
\end{split}
\end{equation}
the points $K_1$ and $K_2$ on the circle $c_B$ are not in the plane $(x,y)$. Hence, the circles $c_A$ and $c_B$, which are circular fibres in the Hopf fibration,  are disjoint. 

As we can rotate any circular fibre of the unit 3-sphere $\mathcal{T}^3$ to the position of $c_A$, it follows that all circular fibers on $\mathcal{T}^3$ are disjoint.
\paragraph{Stereographic projection.} A reasonable choice to create a map of an ordinary \mbox{2-sphere} (e.g., a map of a reference sphere of the Earth) is to take its stereographic projection from the North Pole onto a tangent plane at the South Pole. Despite the distortion of lengths, angles are preserved (the projection is a conformal mapping). As a result, circles on the \mbox{2-sphere} not passing through the North Pole are projected to circles, and circles through the North Pole become straight lines. 
In our analogical four-dimensional case, each point of $\mathcal{T}^3$ is projected via projecting rays through a fixed point  of $\mathcal{T}^3$ -- the center of projection -- onto a \mbox{3-space} touching $\mathcal{T}^3$ at the point antipodal to the center of projection. As this stereographic projection of $\mathcal{T}^3$ to a \mbox{3-space} is conformal, all the circular Hopf fibers are projected to circles apart from the fiber through the center of projection, which projects to a line. We have already shown that the inverse images of points on $\mathcal{B}^2$ in the Hopf fibration are disjoint circles. Moreover, under the stereographic projection, these disjoint fibers are projected onto linked circles (or a line) (see \cite{Lyons2003}). Given a circle on $\mathcal{B}^2$, its inverse image under the Hopf fibration projects to a family of disjoint circles on $\mathcal{T}^3$, and these are projected in the stereographic projection to a torus. Consequently, the stereographic images of all the fibers create nested tori; see \cite{Tsai2004} for Lun-Yi Tsai's exquisite artistic geometric illustrations.

\section{Preliminary constructions}
\label{sec:preliminary}
\subsection{Double orthogonal projection}
The double orthogonal projection of \mbox{4-space} onto two mutually perpendicular \mbox{3-spaces} is a generalization of Monge's projection of an object onto two mutually perpendicular planes (see \cite{Zamboj2018a}). We briefly describe the orthogonal projection of a \mbox{2-sphere} onto a plane in order to aid understanding of the double orthogonal projection of a \mbox{3-sphere} onto two mutually perpendicular \mbox{3-spaces} described subsequently.

\begin{figure}[!htb]
\centering
\includegraphics[width=0.4\linewidth]{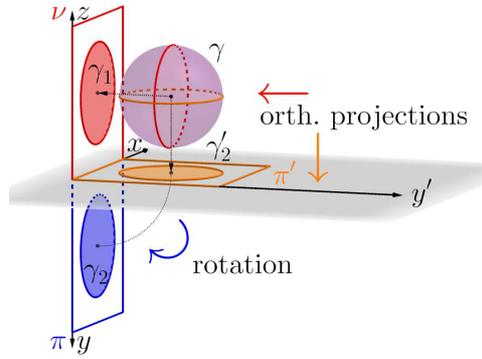}
\caption{Double orthogonal projection of a \mbox{2-sphere} $\gamma$ onto the disk $\gamma_1$ in the 2-space $\nu(x,z)$ and onto the disk $\gamma'_2$ in the 2-space $\pi'(x,y')$ rotated to $\gamma_2$ in $\pi(x,y)$ about the $x$-axis.}
\label{fig:proj3d}
\end{figure}

In an orthogonal projection, the contour generator of a \mbox{2-sphere} is a great circle --- the intersection of the polar plane of the infinite viewpoint with respect to the \mbox{2-sphere} (i.e., the plane perpendicular to the direction of the projection through the center of the \mbox{2-sphere}) and the \mbox{2-sphere} itself. The apparent contour of the \mbox{2-sphere} is also a circle --- the orthogonal projection of the contour generator onto the plane of projection (Fig.~\ref{fig:proj3d}). Therefore, in this three-dimensional case of Monge's projection, we project a \mbox{2-sphere} $\gamma$ to disks $\gamma_1$ and $\gamma'_2$ in two perpendicular planes $(x,z)$ and $(x,y')$. Then the plane $(x,y')$ is rotated about the $x$-axis into the plane $(x,z)$ (the drawing plane) such that the $y'$-axis is rotated to form a $y$-axis that coincides with the $z$-axis but with opposite orientation.

Any point $P$ in the 3-space $(x,y,z)$ is projected orthogonally via its projecting rays to the conjugated image points $P_1$ (front view) and $P_2$ (top view after the rotation). 
The images of projecting rays of $P$ coincide in a single line through $P_1$ and $P_2$, called the ordinal line (or line of recall) of $P$, and it is perpendicular to the $x$-axis in the drawing plane.

\begin{figure}[!htb]
\centering
\includegraphics[width=0.4\linewidth]{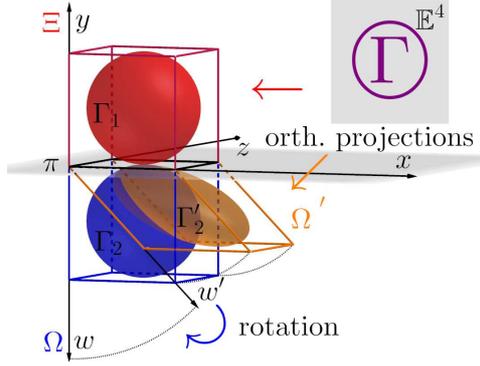}
\caption{Apparent contours of a \mbox{3-sphere} $\Gamma$ in the double orthogonal projection. The \mbox{2-sphere} $\Gamma_1$ in the \mbox{3-space} $\Xi(x,y,z)$ and the \mbox{2-sphere} $\Gamma'_2$ in the \mbox{3-space} $\Omega'(x,z,w')$ rotated to $\Gamma_2$ in $\Omega(x,z,w)$ about the plane $\pi(x,z)$.}
\label{fig:proj4d}
\end{figure}

Analogously, we project a \mbox{3-sphere} in \mbox{4-space} (Fig.~\ref{fig:proj4d}) orthogonally onto a \mbox{3-space}. The contour generator of the \mbox{3-sphere} is a \mbox{2-sphere} -- the intersection of the polar space of the infinite viewpoint with respect to the \mbox{3-sphere} (i.e. the \mbox{3-space} perpendicular to the direction of the projection through the center of the \mbox{3-sphere}) and the \mbox{3-sphere} itself. The apparent contour of the \mbox{3-sphere} is a \mbox{2-sphere} -- the orthogonal projection of the contour generator (\mbox{2-sphere}) onto the \mbox{3-space} of projection. Thus, in the double orthogonal projection both apparent contours in the \mbox{3-spaces} $\Xi(x,y,z)$ and $\Omega'(x,z,w')$ of the \mbox{3-sphere} $\Gamma$ are \mbox{2-spheres} $\Gamma_1$ and $\Gamma'_2$, respectively. Then we rotate $\Omega'(x,z,w')$ about the plane $\pi(x,z)$ to the \mbox{3-space} $\Xi(x,y,z)$ (the modeling \mbox{3-space}) such that the $w'$-axis is rotated to form a $w$-axis that coincides with the $y$-axis but with opposite orientation. After the rotation, we say that $\Gamma$ has two conjugated images in the modeling \mbox{3-space} -- the $\Xi$-image $\Gamma_1$ and the $\Omega$-image $\Gamma_2$ (rotated $\Gamma'_2$).
	
\subsection{Visualization of a point in $\mathbb{R}^4$ and $\mathbb{C}^2$}
\begin{figure}[!htb]
\centering
\includegraphics[width=0.3\linewidth]{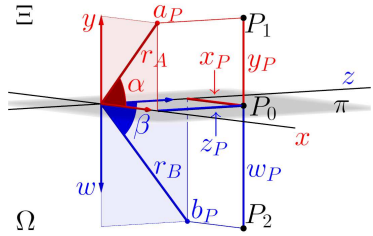}
\caption{Conjugated images $P_1,P_2$ of a point $P(x_P,y_P,z_P,w_P)$ in the double orthogonal projection, and visualization of the complex coordinates of $P$.}
\label{fig:point}
\end{figure}

A point in \mbox{4-space} is given by its two conjugated images ($\Xi$ and $\Omega$-image) lying on their common ordinal line, i.e. the line of coinciding rays of projection after the rotation in the modeling \mbox{3-space} perpendicular to $\pi(x,z)$. In the other direction, let $P$ be a point with coordinates $[x_P,y_P,z_P,w_P]$ in $\mathbb{R}^4$. The $\Xi$-image $P_1[x_P,y_P,z_P]$ and $\Omega$-image $P_2[x_P,z_P,w_P]$ are synthetically constructed in Fig.~\ref{fig:point} using true lengths on the coordinate axes.\footnote{Concretely, in the modeling \mbox{3-space} with the orthogonal coordinate system $(x,y,z)$ given in Fig.~\ref{fig:point}, we should say that the images of $P$ have coordinates $P_1[x_P,y_P,z_P]$ and $P_2[x_P,-w_P,z_P]$. In an implementation such coordinates would naturally be used if we wanted to construct the images from the parametric representation.} The ordinal line of the point $P$ is perpendicular to $\pi(x,z)$ through $P_1$ and $P_2$. Moreover, if we interpret the point $P$ to be in $\mathbb{C}^2$ with coordinates $P[a_P,b_P]$, then $a_P=[x_P,y_P]$ is on the complex line generated by the real axes $x,y$, and $b_P=[z_P,w_P]$ is on another complex line generated by $z,w$. Using the trigonometric representation, we have $a_P=[r_A\cos\alpha,r_A\sin\alpha]$ and $b_P=[r_B\cos\beta,r_B\sin\beta]$ (cf. equation~(\ref{eq:parS3gon})).

\subsection{Visualization of a point on a \mbox{3-sphere}}
\begin{figure}[!htb]
\centering
\subfloat[]{\includegraphics[height=0.4\textheight]{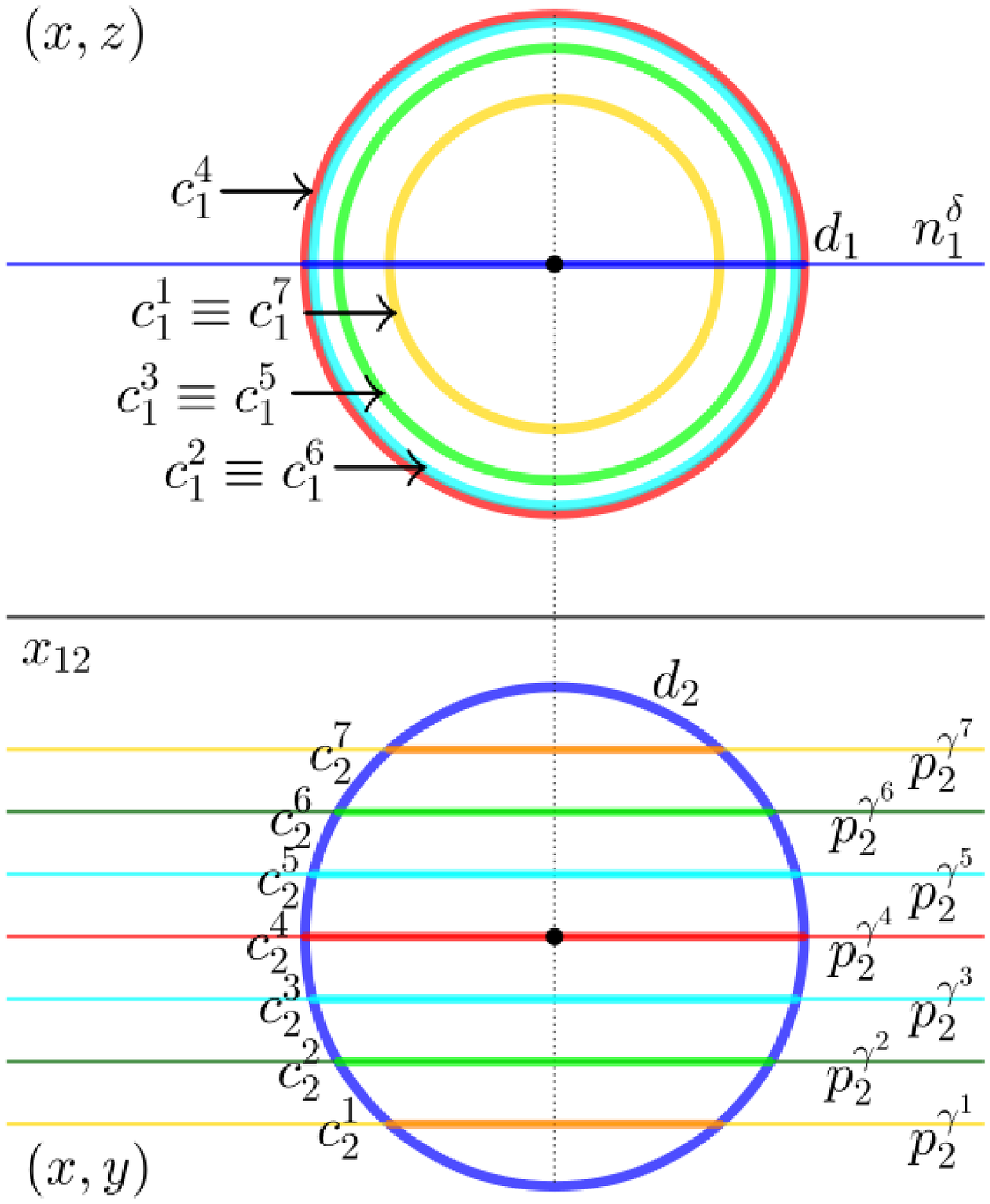}%
\label{fig:circle_sections}}
\hfil
\subfloat[]{\includegraphics[height=0.4\textheight]{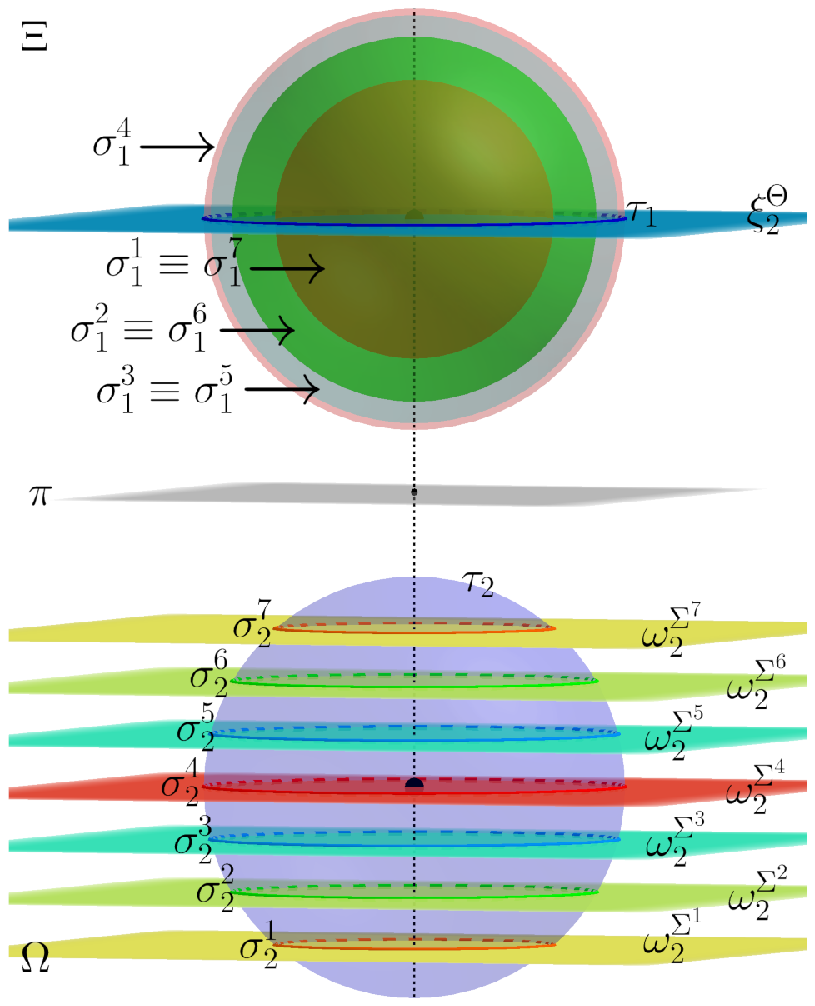}%
\label{fig:2sphere_sections}}
\caption{(a) Circular sections $c^1,\dots,c^7$ of a 2-sphere with planes $\gamma^1,\dots,\gamma^7$ parallel to $(x,z)$ in Monge's projection. The planes are given by their intersections $p^\gamma$ with the plane $(x,y)$. The apparent contours of the sections are circles $c_1^1,\dots,c_1^7$ in the front view and segments $c_2^1,\dots,c_2^7$ in the top view. (b) Sections $\sigma^1,\dots,\sigma^7$ of a 3-sphere with 3-spaces $\Sigma^1,\dots,\Sigma^7$ parallel to the 3-space $\Xi(x,y,z)$ are 2-spheres. The 3-spaces are given by their intersections $\omega^\Sigma$ 
with the 3-space $\Omega$. The $\Xi$-images $\sigma_1^1,\dots,\sigma_1^7$ of the spherical sections are 2-spheres in their true shape and $\Omega$-images $\sigma_2^1,\dots,\sigma_2^7$ are disks.}
\label{fig:sphere_sections}
\end{figure}

In an orthogonal projection of a 2-sphere onto a plane, points on the 2-sphere are projected onto a disk. To locate the position of a point in the image, we can draw its circle of latitude (in a plane perpendicular to the direction of the projection). Analogously, when we orthogonally project a \mbox{3-sphere} onto a \mbox{3-space}, we can locate a point on the \mbox{3-sphere} by drawing its ``\mbox{2-sphere} of latitude'' (in a \mbox{3-space} perpendicular to the direction of the projection). In Monge's projection, planes parallel to the plane $(x,z)$ intersect the plane $(x,y)$ in lines parallel to $x$. The sections of a 2-sphere with planes parallel to $(x,z)$ in Monge's projection are shown in Fig.~\ref{fig:circle_sections}; they create line segments in the top view and circles at their true size in the front view. In other words, if we imagine a 2-sphere passing orthogonally through a plane, their common intersection will be the tangent point extending continuously to the great circle and shrinking back to a point. Analogously, the 3-spaces parallel to the 3-space $\Xi(x,y,z)$ intersect the 3-space $\Omega(x,z,w)$ in planes parallel to $\pi(x,z)$. If a 3-sphere passes orthogonally through a 3-space, their common intersection is the tangent point, extending to the great 2-sphere and shrinking back to a point. Therefore, the intersections of the 3-sphere with a bundle of 3-spaces parallel to $\Xi(x,y,z)$ are 2-spheres (Fig.~\ref{fig:2sphere_sections}). Their $\Xi$-images are circles and $\Omega$-images are 2-spheres at their true size. This construction should give us some insight into the visualization of points on tori inside a 3-sphere that will be carried out later. 

\begin{figure}[!htb]
\centering
\includegraphics[width=0.4\linewidth]{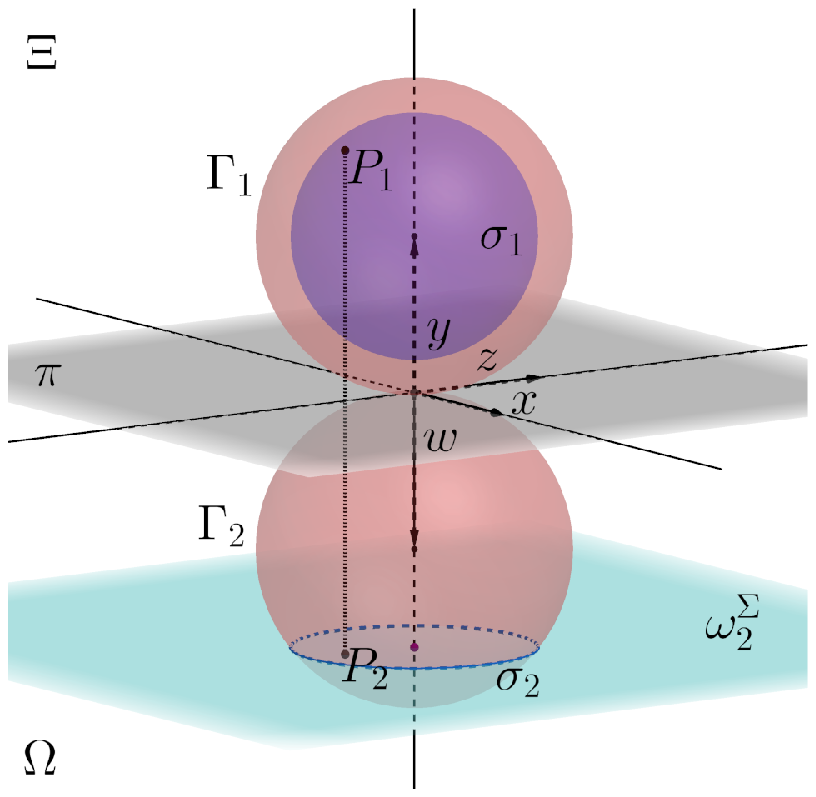}
\caption{Conjugated images $P_1$ and $P_2$ of a point $P$ on the spherical section $\sigma$ of a \mbox{3-sphere} $\Gamma$ with a \mbox{3-space} $\Sigma$ parallel to $\Xi(x,y,z)$. The $\Xi$-image $P_1$ is on $\sigma_1$, which is in its true shape. The $\Omega$-image $\sigma_2$ of the 2-sphere $\sigma$ is a disk with the same radius as $\sigma_1$. The $\Omega$-image $P_2$ lies on the disk $\sigma_2$.
\newline
In the interactive 3D model \texttt{https://www.geogebra.org/m/yt27evc8} (or Suppl. File~1), the user can manipulate the position of $P_1$ on $\sigma_1$ and $\omega_2^\Sigma$ with the ``Moving point'' on the $w$-axis and observe the conjugated images of positions of the point $P$ on the 3-sphere.}
\label{fig:pointionsphere}
\end{figure}

Let $P$ be a point on a \mbox{3-sphere} $\Gamma$ (Fig.~\ref{fig:pointionsphere}, Suppl. File~1). The point $P$ lies on some \mbox{2-sphere} $\sigma$ in a \mbox{3-space} $\Sigma$ parallel to $\Xi(x,y,z)$. Since $\Sigma$ has the same $w$-coordinate as $P$, its $\Omega$-image appears as the plane $\omega^\Sigma_2$ through~$P_2$ parallel to $\pi(x,z)$. The intersection of $\omega^\Sigma_2$ and $\Gamma_2$ is a circle~$\sigma_2$ that is the boundary of the $\Omega$-image of $\sigma$. The $\Xi$-image $\sigma_1$ through~$P_1$ is a \mbox{2-sphere} concentric with $\Gamma_1$ and has radius equal to the radius of $\sigma_2$.

\subsection{Stereographic projection of a point on a \mbox{3-sphere} onto a \mbox{3-space}}
\label{sec:stproj}
\begin{figure}[!htb]
\centering
\includegraphics[width=0.6\linewidth]{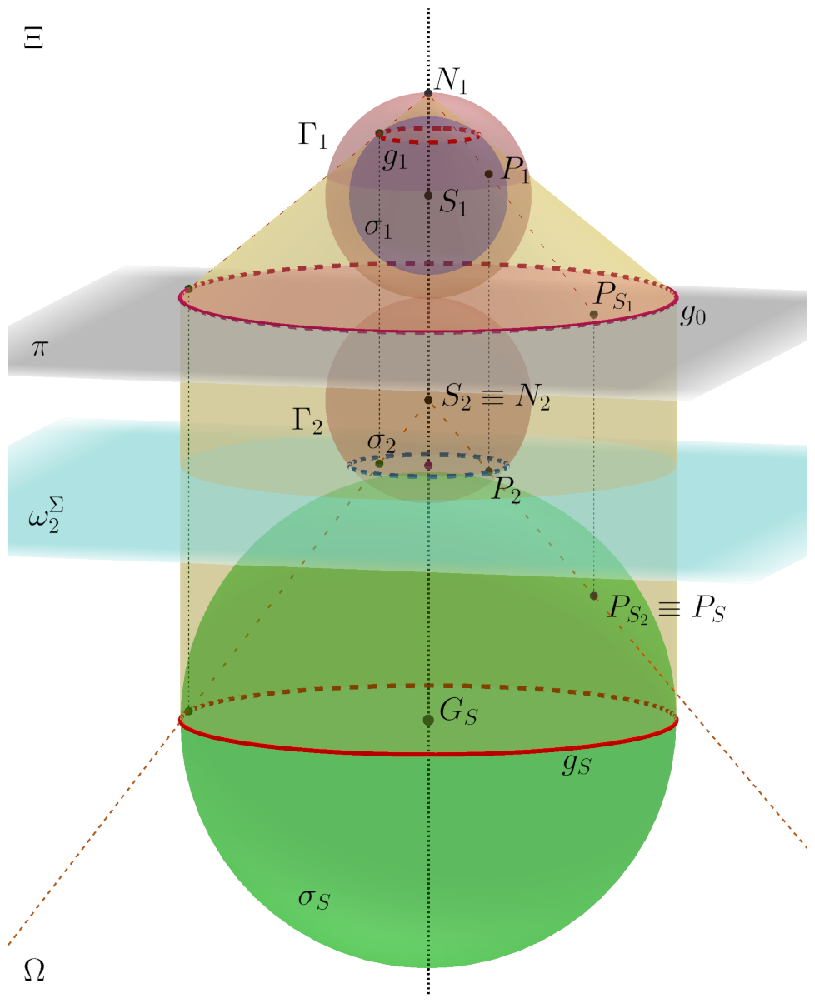}
\caption{The stereographic image $P_S$ of a point $P$ on a \mbox{3-sphere} $\Gamma$ projected from the center $N$ to the tangent \mbox{3-space} $\Omega(x,z,w)$.  The stereographic image $P_S$ lies in the 3-space $\Omega$, so it is also the $\Omega$-image $P_{S_2}$ constructed as the intersection of $NP$ and $\Omega$. Additionally, the spherical section $\sigma$ from Fig.~\ref{fig:pointionsphere} is stereographically projected to the 2-sphere $\sigma_S$. For the analytic coordinates of $P_S$ see equation~(\ref{eq:parS3gonSter+}).\newline
In the interactive 3D model \texttt{https://www.geogebra.org/m/xth6uszb} (or Suppl. File~2), the user can, again, manipulate $P_1$ on $\sigma_1$ and $\omega_2^\Sigma$ by ``Moving point'' on $w$-axis and observe the stereographic images.}
\label{fig:sterproj}
\end{figure}

	Let us have a unit \mbox{3-sphere} $\Gamma$ with the center $S=[0,1,0,1]$ (cf. equation~(\ref{eq:3sphere}))
\begin{equation}
	\label{eq:3sphere2}
	x^2+(y-1)^2+z^2+(w-1)^2=1.
	\end{equation}
	We project points of the \mbox{3-sphere} $\Gamma$ from its point $N=[0,2,0,1]$ onto the \mbox{3-space} $\Omega(x,z,w): y=0$ that touches~$\Gamma$ at the antipodal point $M=[0,0,0,1]$. Let $P$ be a point on $\Gamma$ with images $P_1$ and $P_2$ (Fig.~\ref{fig:sterproj}, Suppl. File~2). The stereographic image $P_S$ of the point $P$ in $\Omega(x,z,w)$ is the intersection of the line $NP$ and $\Omega(x,z,w)$. The line $N_1P_1$ intersects $\pi(x,z)$ in the point $P_{S_1}$. The intersection of the ordinal line through $P_{S_1}$ and the $\Omega$-image $N_2P_2$ is the desired point $P_{S_2}$ that coincides with $P_S$ in \mbox{4-space}. 
	
Since stereographic projection is a conformal mapping, we can also conveniently project the \mbox{2-sphere} $\sigma$ through the point $P$ described in the previous section. The stereographic image $\sigma_S$ of the \mbox{2-sphere} $\sigma$ is a \mbox{2-sphere} $\sigma_S$ with center $G_S$ on the ordinal line through the center $S_1$ of $\sigma_1$. Its equatorial circle $g_S$ is the image of the tangent circle $g_1$ on the sphere $\Gamma_1$ of the cone with the vertex $N_1$, i.e. the intersection of the polar plane of the pole $N_1$ with respect to $\sigma_1$ and $\sigma_1$ itself. Therefore, with the use of a point on $g_1$, we construct the image circle $g_S$ of the circle $g$, and the center $G_S$ of $g_S$ is the center of the \mbox{2-sphere} $\sigma_S$ -- the stereographic image of the \mbox{2-sphere} $\sigma$. 

\section{Synthetic construction of a Hopf fiber}
\label{sec:construction}
\begin{figure}[!htb]
\centering
\includegraphics[width=0.6\linewidth]{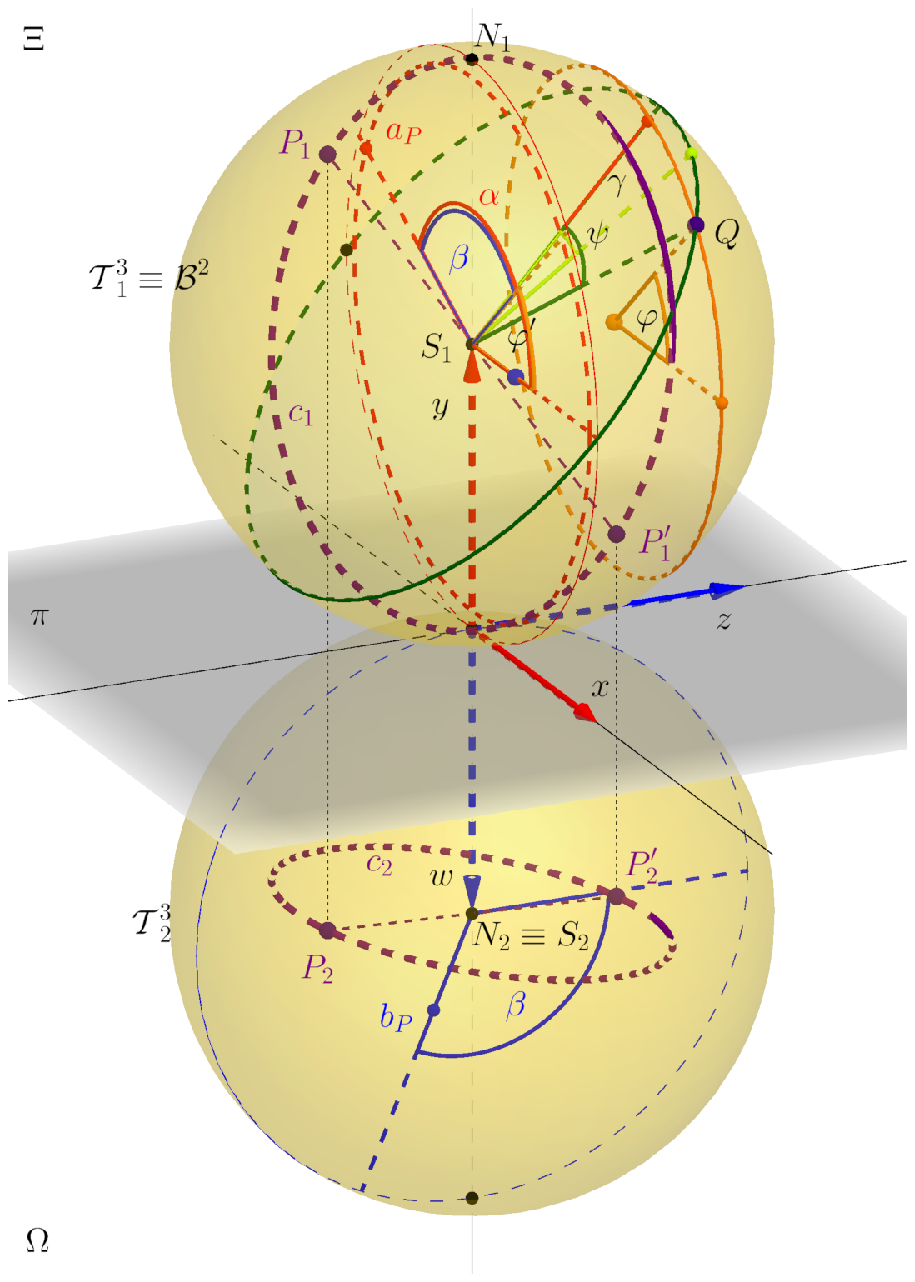}
\caption{Construction of the Hopf fiber $c$ on $\mathcal{T}^3$ corresponding to a point $Q$ on $\mathcal{B}^2$.\newline
With the step-by-step construction \texttt{https://www.geogebra.org/m/w2kugajz} (or Suppl. File~3), the reader can follow the steps of the construction from the text in the {\sl Construction Protocol} window. Use the arrows to move between the steps. From step 1 onward, the user can move the point $Q$ on $\mathcal{B}^2$ and dynamically change all the dependent elements of the model and the resulting fiber. From step 3 onward, the user can also manipulate the angle $\beta$ and observe positions of the conjugated images of $P$ and $P'$ on the fiber. Steps 10 and 11 give the construction of the stereographic projection in Fig.~\ref{fig:hopf_SBS}.}
\label{fig:hopf_SBS1}
\end{figure}

\begin{figure}[!htb]
\centering
\subfloat[]{\includegraphics[width=0.16\textwidth]{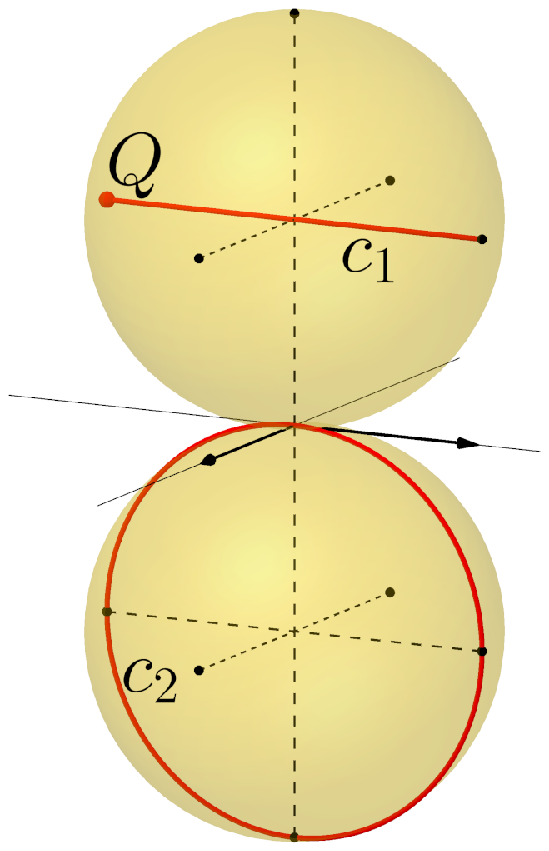}%
\label{fig:hopf_fiber1}}
\subfloat[]{\includegraphics[width=0.16\textwidth]{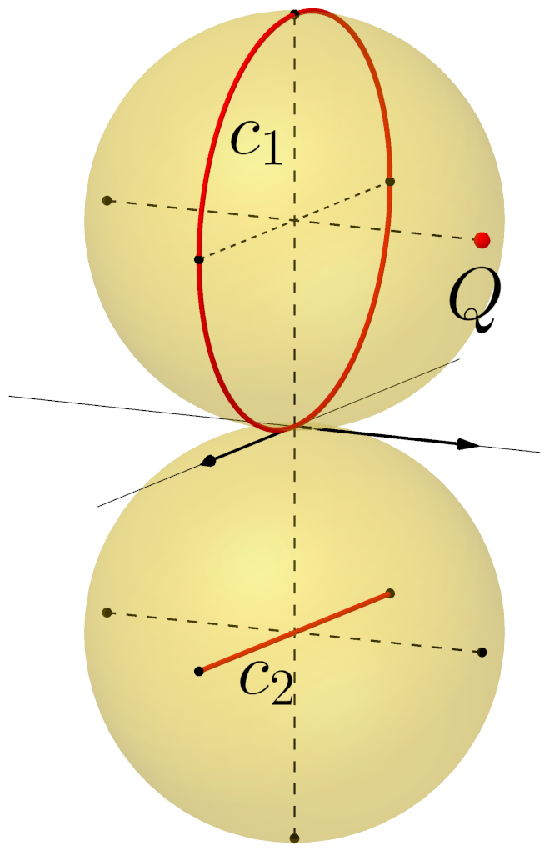}%
\label{fig:hopf_fiber2}}
\subfloat[]{\includegraphics[width=0.16\textwidth]{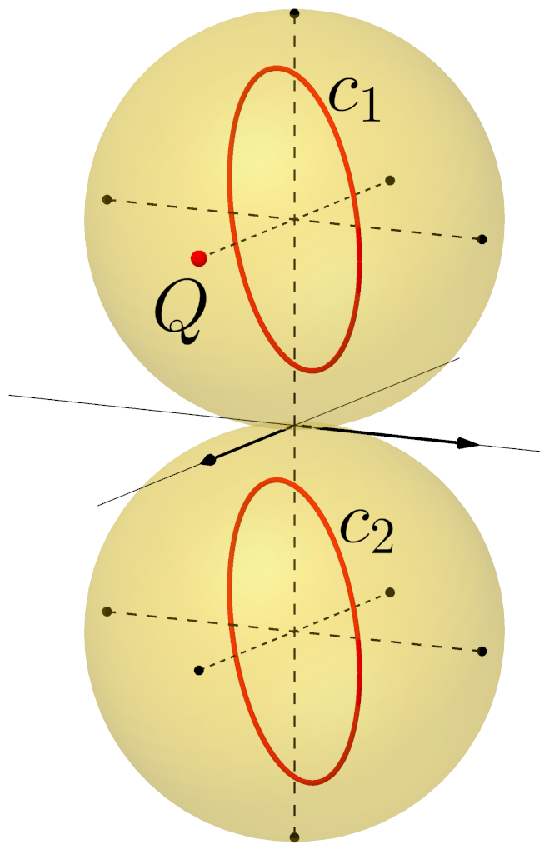}%
\label{fig:hopf_fiber3}}
\subfloat[]{\includegraphics[width=0.16\textwidth]{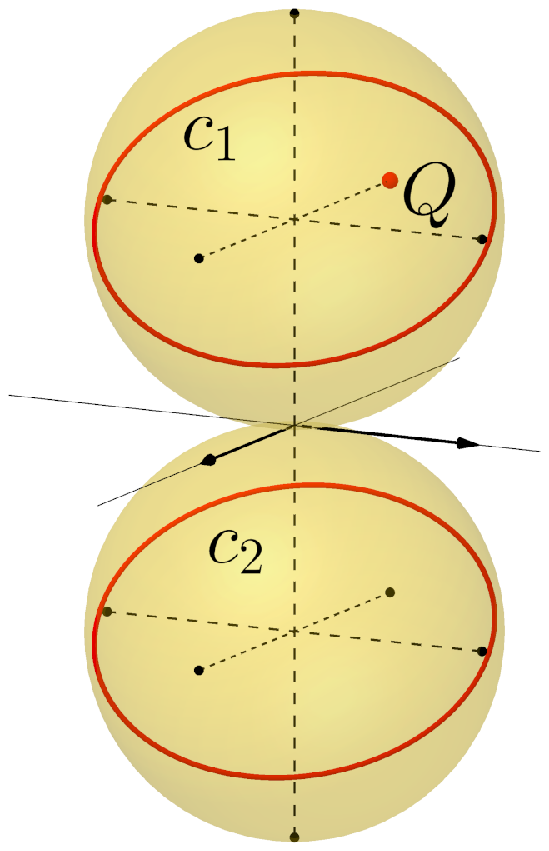}%
\label{fig:hopf_fiber5}}
\subfloat[]{\includegraphics[width=0.16\textwidth]{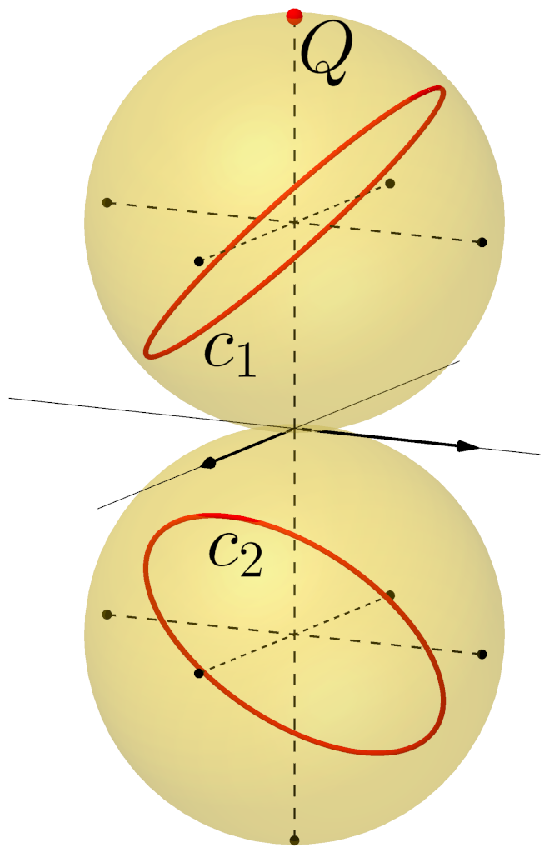}%
\label{fig:hopf_fiber4}}
\subfloat[]{\includegraphics[width=0.16\textwidth]{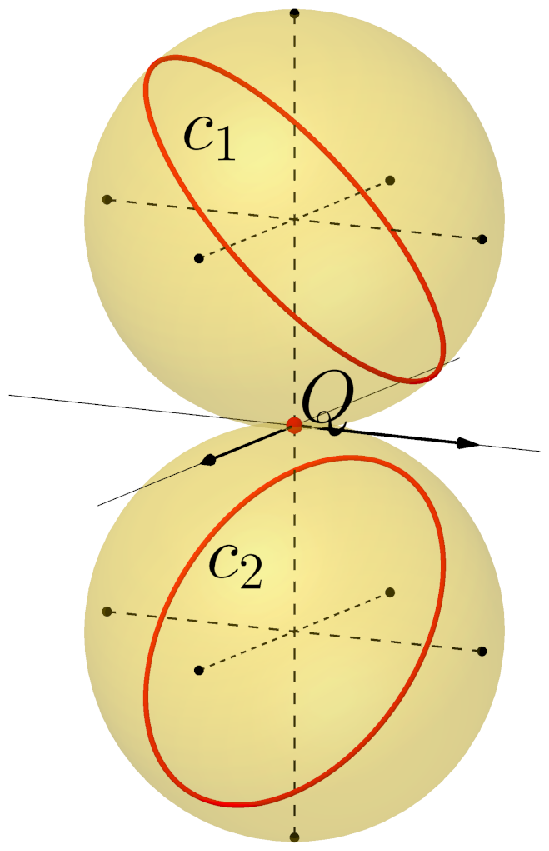}%
\label{fig:hopf_fiber6}}
\caption{Conjugated images of the fibers corresponding to the point Q in certain positions: (a) Q=[0,1,-1], (b) Q=[0,1,1], (c) Q=[1,1,0], (d) Q=[-1,1,0], (e) Q=[0,2,0], (f) Q=[0,0,0].}
\label{fig:hopf_fiber_positions}
\end{figure}

\begin{figure*}[!htb]
\centering
\subfloat[]{\includegraphics[width=0.3\linewidth]{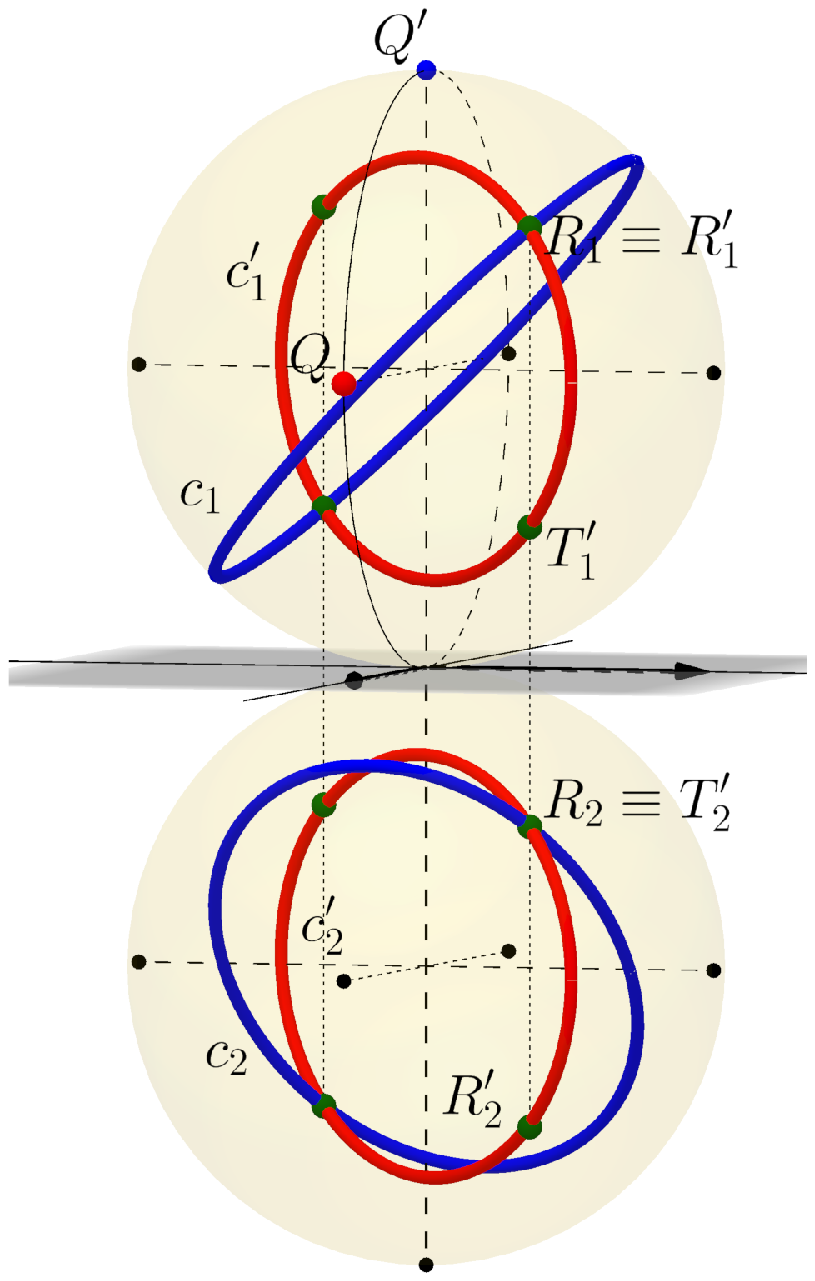}
\label{fig:hopf_special}}
\hfil
\subfloat[]{\includegraphics[width=0.3\linewidth]{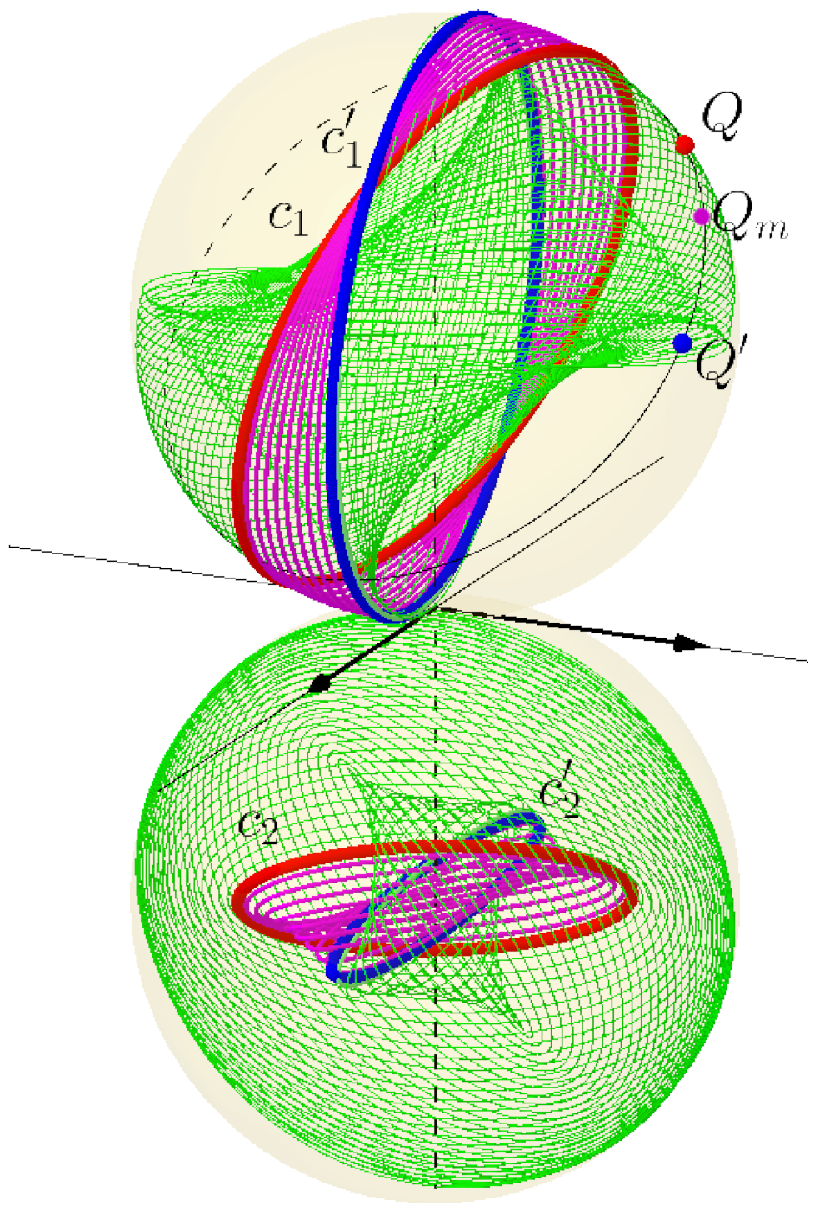}
\label{fig:hopf_trace}}
\caption{(a) Special position of the points $Q$ and $Q'$ such that the conjugated images of disjoint fibers $c$ and $c'$ intersect. In the $\Xi$-image, $c_1$ and $c'_1$ intersect in a point $R_1\equiv R'_1$, but their $\Omega$-images $R_2\in c_1$ and $R'_2\in c'_2$ are distinct. Similarly the point $T'$ on $c'$ has its $\Omega$ image $T'_2\equiv R_2$, but their $\Xi$-images are distinct. (b) Motion of the point $Q_m$ between $Q$ and $Q'$. The pink fibers correspond to positions of $Q_m$ on the shorter arc between $Q$ and $Q'$, and the green fibers correspond to the motion of $Q_m$ on the longer arc.\newline
In the interactive model \texttt{https://www.geogebra.org/m/ebjkx8pj} (or Suppl. File~4), the user can with the scroll bars dynamically change the positions of $Q$ and $Q'$ dependent (via equation~(\ref{eq:parS2+})) on parameters $(\varphi, \psi)$ and $(\varphi', \psi')$, respectively. The point $Q_m$ can move freely along the great circle through $Q$ and $Q'$. The conjugated images of the fibers corresponding to $Q_m$ draw their traces in the modeling 3-space.}
\label{fig:hopf2fibers}
\end{figure*}

In this section, we 
illustrate geometrically the mathematical properties of the Hopf fibration given in Section~\ref{sec:math} in images under the double orthogonal projection. This way, we synthetically construct a circular fiber on a \mbox{3-sphere} from a point on a \mbox{2-sphere} by elementary geometric constructions, thereby liberating it from its analytic description. For the puposes of computation we used the equations of a 3-sphere $\mathcal{T}^3$ with the center at the origin, but for the sake of visualization (to differentiate the $\Xi$ and $\Omega$-images) it proves more suitable to use a 3-sphere $\mathcal{T}^3$ with center $[0,1,0,1]$. Then the \mbox{2-sphere} $\mathcal{B}^2$ has center with coordinates $[0,1,0]$ in $\Xi(x, y, z)$. Such a translation does not influence the properties of the Hopf fibration. This applies for Figs.~\ref{fig:hopf_SBS1}--\ref{fig:hopf2fibers}, \ref{fig:hopf_SBS}--\ref{fig:nested_tori_S}; the parametric equations corresponding to the visualizations are in the appendix.

Let $Q\in\Xi(x,y,z)$ be an arbitrary point on $\mathcal{B}^2$ (Fig.~\ref{fig:hopf_SBS1}, Suppl. File~3). We will find its Hopf fiber -- a circle~$c$ on $\mathcal{T}^3$. 
Equations~(\ref{eq:spherical}), \ref{eq:parS2+} and (\ref{eq:hopfcoordinates},\ref{eq:parS3gon+}) give a relation between the spherical coordinates of the point $Q$ (with parameters $\varphi$ and $\psi$) and the Hopf coordinates of the fiber $c$ (with parameters $\varphi,\psi$ and $\beta$). The construction proceeds by the following steps (see the step-by-step construction \texttt{https://www.geogebra.org/m/w2kugajz} (or Suppl. File~3)):
\begin{enumerate}
\newcounter{steps}
	\item Construct any point $Q$ on $\mathcal{B}^2$.
	\item Find the angle $\varphi$ (from equation~(\ref{eq:spherical})): Construct the plane parallel to $(x,y)$ through $Q$ that cuts $\mathcal{B}^2$ in a circle. The oriented angle between the radius parallel to the \mbox{$x$-axis} and the radius terminating in the point $Q$ is the angle $\varphi$.
	\item Construct an arbitrary angle $\beta$ such that we can graphically add it to $\varphi$: First, translate $\varphi$ to $\varphi'$ in the plane $(x,y)$ with its vertex in the center of $\mathcal{B}^2$, and the initial side in the direction of the non-negative \mbox{$x$-axis}. Now choose $\beta$ with the same vertex and initial side in the terminal side of $\varphi'$. For implementation, it is enough to choose $\beta$ on the top semicircle, as the points $P$ and $P'$ constructed in step 7 dependent on $\beta$ will be antipodal. Additionally, construct $\alpha$ such that $\alpha=\varphi'+\beta$ (cf. equation~(\ref{eq:spherical})). The angles $\alpha$ and $\beta$ are arguments of the complex points $a_P=[r_A\cos\alpha,r_A\sin\alpha]\in (x,y)$ and $b_P=[r_B\cos\beta,r_B\sin\beta]\in(z,w)$ (see the derivation of $z_1$ and $z_2$ above equation~(\ref{eq:parS3gon})).
	\item Find the angle $\psi$ (from equation~(\ref{eq:spherical})): Construct the great circle of $\mathcal{B}^2$ through $Q$ with a diameter parallel to the $x$-axis. Choose a radius in the $(y,z)$-plane to be the initial side of the angle $\psi$ with the terminal side being the radius through $Q$.
	\item Construct the moduli of the points $a_P$ and $b_P$ (cf. equation~(\ref{eq:parS3gon})): Let $\gamma$ be the half-angle of $\psi$ and find its cosine by dropping a perpendicular onto the initial side of $\psi$. The length $\cos\frac{\psi}{2}=\cos\gamma$ is the modulus $r_A$ of the point $a_P$. Similarly find the length $\sin\frac{\psi}{2}=\sin\gamma$ on the radius in the direction of the $x$-axis, which is the modulus $r_B$ of $b_P$.
	\item Construct points $a_P$ and $b_P$: Using the moduli and arguments of $a_P$ and $b_P$, construct them according to Fig.~\ref{fig:point}.
\item Construct the $\Xi$ and $\Omega$-images of the point $P$: Having $a_P$ and $b_P$, we finalize the images $P_1$ and $P_2$ on the parallels to the reference axes (as in Fig.~\ref{fig:point}).
\item Construct the antipodal point $P'$ on $\mathcal{B}^2$: Parallel projection preserves central symmetry, and the images $P'_1$ and $P'_2$ are the reflections of $P_1$ and $P_2$ about the centers $S_1$ and $S_2$ of the images of $\mathcal{T}^3$.
\item Construct the Hopf fiber corresponding to the point $Q$: The Hopf fiber is a circle $c$ consisting of the locus of points $P$ dependent on $\beta$. We have a point construction of $P$, which can be repeated for different choices of $\beta$ even though, at this point, we comfortably use the GeoGebra tool to draw the locus. The orthogonal projections of $c$ will appear as ellipses (or circles, or line segments) $c_1$ and $c_2$.
	\setcounter{steps}{\value{enumi}}
\end{enumerate}
Varying $\beta$ in the interactive applet, $P$ moves on its fiber $c$. Moving with $Q\in\mathcal{B}^2$, the whole fiber $c$ moves on $\mathcal{T}^3$. 

To illustrate the construction, we consider the images of fibers corresponding to certain points on the base sphere $\mathcal{B}^2$ on diameters parallel to the $x,y,$ and $z$ axes (Fig. \ref{fig:hopf_fiber_positions}). 
\begin{itemize}
\item[] (with respect to the translated coordinates of the centers of $\mathcal{B}^2$ and $\mathcal{T}^3$)  
\item[] (a), (b) The fibers of points $[0,1,\pm1]$ (i.e. $\varphi=0$ and $\psi=0,\pi$) are in the planes $(x,w)$ and $(y,z)$, respectively. Thus, their conjugated images are a line segment and a great circle. In particular, for $\psi=\pi$ the point on the base sphere lies on its fiber.
\item[] (c), (d) The fibers of points $[\pm1,1,0]$ (i.e. $\varphi=0, \pi$ and $\psi=\frac{\pi}{2}$) lie in the plane of symmetry of the $x$ and $z$ axes and their conjugated images are congruent.
\item[] (e), (f) The fibers of points $[0,1+\pm1,0]$ (i.e. $\varphi=\frac{\pi}{2}, \frac{3\pi}{2}$ and $\psi=\frac{\pi}{2}$) have their $\Xi$-images in the plane of symmetry of the $y$ and $z$ axes and $\Omega$-images in the plane of symmetry of the $x$ and $w$ axes.
\end{itemize}

\begin{figure}[!htb]
\centering
\includegraphics[width=0.6\linewidth]{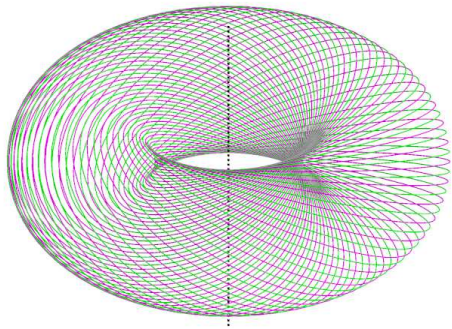}
\caption{A torus in a three-dimensional space generated by revolution of Villarceau circles. All the generating circles are interlinked.}
\label{fig:torusflips}
\end{figure}

We have already mentioned that fibers corresponding to two distinct points on the base 2-sphere $\mathcal{B}^2$ create linked circles. Let $Q$ and $Q'$ be two points on $\mathcal{B}^2$ and $c$ and $c'$ their fibers, respectively. We can easily observe from the conjugated images of $c$ and $c'$ that the fibers are disjoint, for if $c$ and $c'$ had a point of intersection $R$, their conjugated images $c_1, c'_1$ and $c_2,c'_2$ would intersect in the conjugated images $R_1, R_2$ of the point $R$. In the case that the conjugated images of a fiber are in a plane perpendicular to $\pi$ (e.g., Fig.~\ref{fig:hopf_special}), we must not swap the conjugated images of its points. Furthermore, we should note that the circles are linked on the 3-sphere $\mathcal{T}^3$, but this property cannot be validated in the embedding 4-space. Analogously, imagine a point in a circular region on a 2-sphere embedded in 3-space. On the 2-sphere, the point cannot escape the bounding circle. However, if we remove the 2-sphere and leave only the 3-space, the point is not bounded at all. Hence, to establish the fibers' interlinkedness, we need to understand the topology of the underlying 3-sphere $\mathcal{T}^3$. Let us construct a great circle on $\mathcal{B}^2$ through $Q$ and $Q'$ and observe the motion of the fiber of a point $Q_m$ on the circle moving from $Q$ to $Q'$ (Fig.~(\ref{fig:hopf_trace})). During this motion, the moving fiber twists along a surface. In fact, the generating fibers are always one of a pair of Villarceau circles around a torus (cf. Fig.~\ref{fig:torusflips} with a 3-dimensional parallel projection onto a plane of a torus generated by Villarceau circles). The toroidal structure of a 3-sphere will be seen more clearly in Section~\ref{sec:tori} with the use of stereographic projection.

\subsection{Construction of the stereographic image of a Hopf fiber}

\begin{figure}[!htb]
\centering
\includegraphics[width=0.6\linewidth]{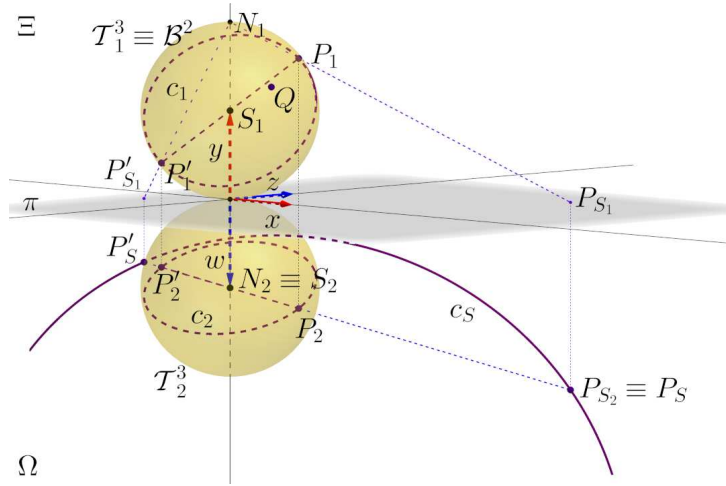}
\caption{Construction of the stereographic image $c_S$ of the Hopf fiber $c$ (cf. Fig.~\ref{fig:sterproj}). The fiber $c_S$ is constructed as the locus of points $P_S$ and $P'_S$ dependent on the angle $\beta$.\newline
Stereographic images are constructed in steps 10 and 11 in the step-by-step construction \texttt{https://www.geogebra.org/m/w2kugajz} (or Suppl. File~3) described in Fig.~\ref{fig:hopf_SBS1}.}
\label{fig:hopf_SBS}
\end{figure}

To see the circular structure of the Hopf fibration, we construct the fibers in stereographic projection (Fig.~\ref{fig:hopf_SBS}). We use the same center of projection and antipodal tangent space $\Omega(x,z,w)$ as in Section~\ref{sec:stproj}. Continuing from the previous construction:
\begin{enumerate}
\setcounter{enumi}{\value{steps}}
\item Construct a stereographic image of the point $P$: Let $N=[0,2,0,1]$ be on $\mathcal{T}^3$ and let the \mbox{3-space} $\Omega(x, z, w)$ be tangent to $\mathcal{T}^3$ at the point $[0,0,0,1]$. The intersection of $N_1P_1$ with $\pi(x,z)$ is $P_{S_1}$. Dropping a perpendicular from $P_{S_1}$ to the line $N_2P_2$ gives us the point $P_{S_2}$ that is also the true stereographic image $P_S$.
\item Construct a stereographic image of the fiber $c$: We use the locus tool from GeoGebra to construct a locus of points $P_S$ dependent on the angle $\beta$. The stereographic image $c_S$ of $c$ is a circle or a line segment (if $c$ passes through $N$). 
\end{enumerate}
Again, by varying $\beta$ the point $P_S$ moves on $c_S$, and by varying the position of $Q\in\mathcal{B}^2$ the circle $c_s$ changes in such a way that it may cover the whole modeling \mbox{3-space}. In the following section, we show how to move $Q$ so as to obtain the Hopf tori. 
\begin{figure}[!htb]
\centering
\subfloat[]{\includegraphics[height=0.25\textheight]{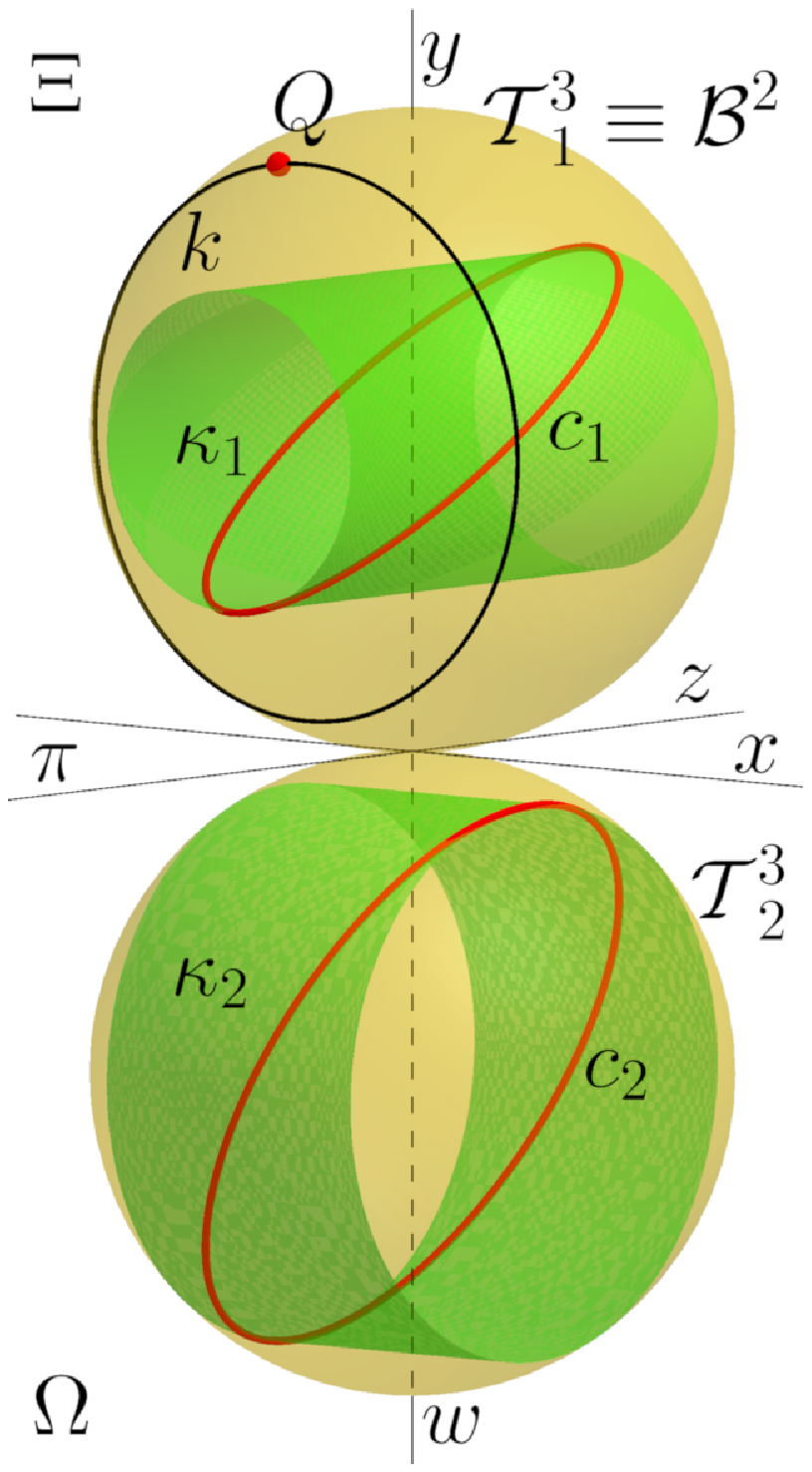}%
\label{fig:hopf_tori1}}
\hfil
\subfloat[]{\includegraphics[height=0.25\textheight]{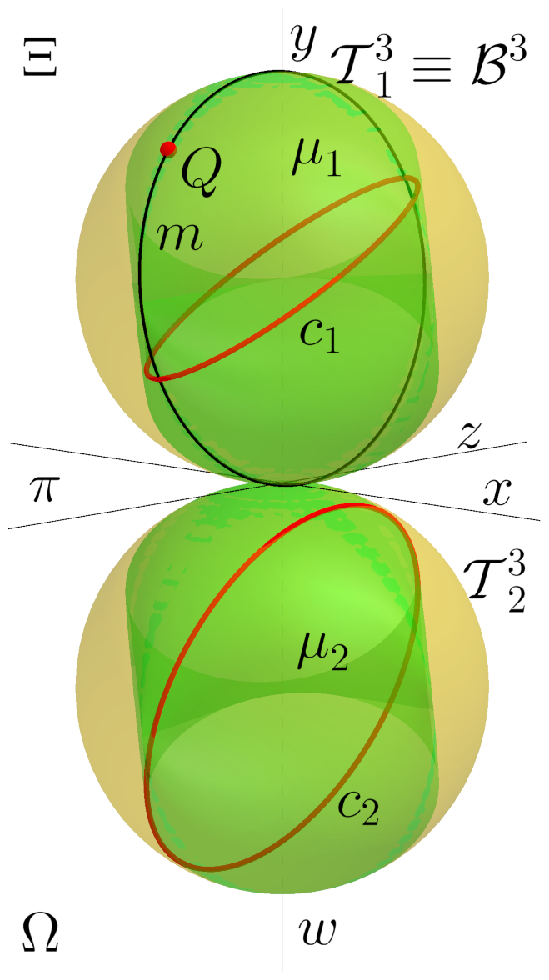}%
\label{fig:hopf_tori2}}
\caption{(a) Torus $\kappa$ on $\mathcal{T}^3$ corresponding to a circle $k$ on $\mathcal{B}^2$ parallel to the $(x,y)$-plane. The torus is generated by fibers $c$ above points $Q$ along the circle $k$. (b) Torus $\mu$ on $\mathcal{T}^3$ corresponding to a circle $m$ on $\mathcal{B}^2$ with a diameter parallel to $z$-axis. Again, the torus is generated by fibers $c$ above points $Q$ on $m$.}
\label{fig:hopf_tori}
\end{figure}

\section{Hopf tori corresponding to circles on $\mathcal{B}^2$}
\label{sec:tori}
\begin{figure*}[!htb]
\centering
\subfloat[]{\includegraphics[height=0.3\textheight]{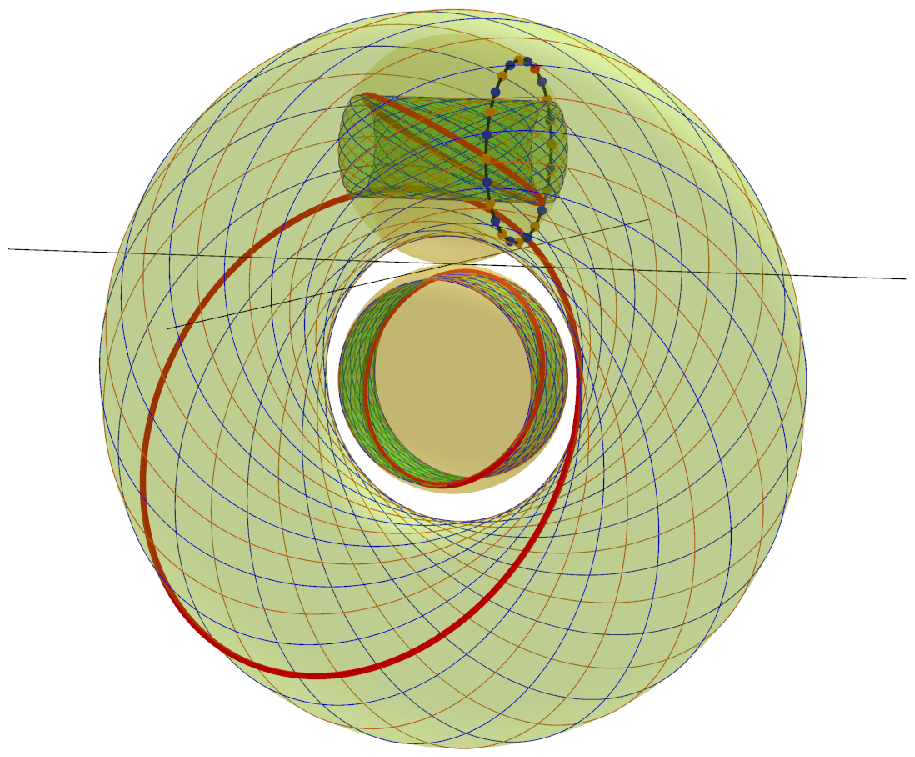}%
\label{fig:hopf_torus1_circ_S}}\hfil
\subfloat[]{\includegraphics[height=0.3\textheight]{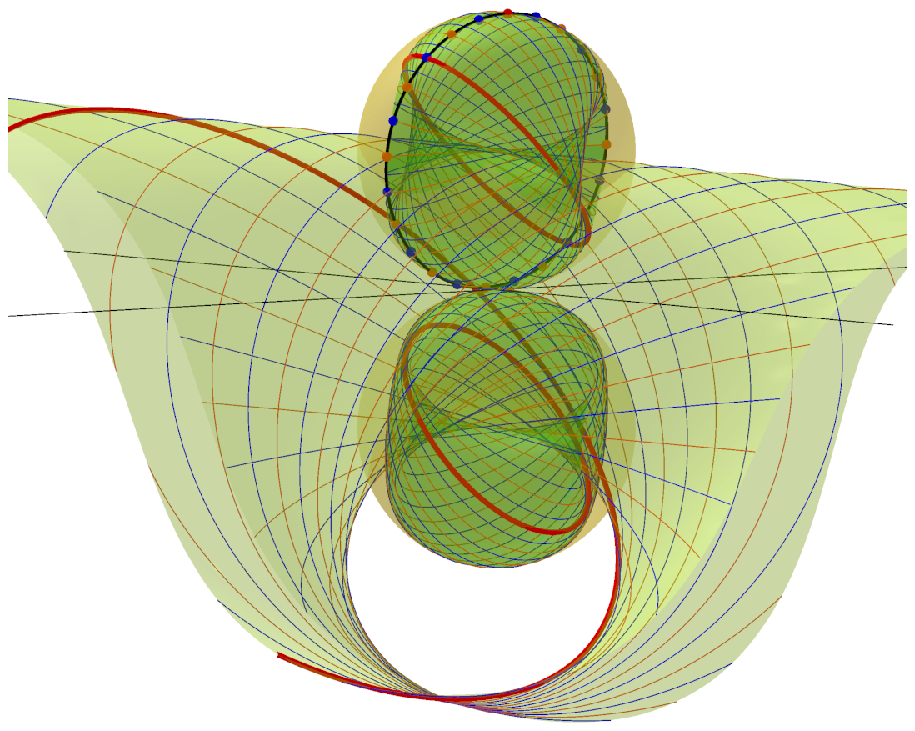}%
\label{fig:hopf_torus2_circ_S}}
\caption{(a) Blue and orange points on a circle on $\mathcal{B}^2$ parallel to the $(x,y)$-plane and the corresponding family of circular fibers on $\mathcal{T}^3$ generating a torus. The stereographic projection of this torus onto $\Omega(x,z,w)$ is a torus of revolution. (b) Blue and orange points on a circle on $\mathcal{B}^2$ with a diameter parallel to the $z$-axis and their fibers. One of the fibers passes through the center of stereographic projection, so its image is a line.
\newline
In the interactive models (a) \texttt{https://www.geogebra.org/m/n4xg3sw6} (or Suppl. File~5) and (b) \texttt{https://www.geogebra.org/m/vrasywpt} (or Suppl. File~6), the user can vary the spherical coordinates of $Q$ in equation~(\ref{eq:parS2+}) and interactively change its circle on $\mathcal{B}^2$ and the corresponding torus. The user can also turn off the visibility of the objects on $\mathcal{B}^2$, the conjugated images in the double orthogonal projection, or the stereographic images.}
\label{fig:hopf_tori_circ_S}
\end{figure*} 
The geometric nature of the Hopf fibration becomes fully apparent when we visualize the tori of Hopf fibers corresponding to circles on the base \mbox{2-sphere} $\mathcal{B}^2$. With respect to the chosen stereographic projection, we divide the following constructions into two cases. First, we construct the tori on $\mathcal{T}^3$ in $\mathbb{R}^4(x,y,z,w)$ corresponding to circles on~$\mathcal{B}^2$ parallel with the plane $(x,y)$ in the \mbox{3-space} $\Xi(x,y,z)$, and then the tori on $\mathcal{T}^3$ corresponding to circles on $\mathcal{B}^2$ with diameter parallel to $z$. Instead of point-by-point constructions, in the following interactive demonstrations (Suppl. Files~5 and~6) the objects are defined by their parametric representations (see equations~(\ref{eq:circ1S2+}) of the circle $k$ on $\mathcal{B}^2$, (\ref{eq:torus1S3+}) of the corresponding torus $\kappa$, (\ref{eq:torus1ster+}) of the stereographic image of the torus $\kappa$ in the appendix), and the user can manipulate the angles $\varphi$ and $\psi$ using the sliders. 
\subsection{Hopf torus of a circle parallel to $(x,y)$}
\label{ss:1}
Let $Q$ be a point on $\mathcal{B}^2$ in $\Xi(x,y,z)$ and $k$ a circle parallel to $(x,y)$ through the point $Q$ (Fig.~\ref{fig:hopf_tori1}, Suppl. File~5). From equation~(\ref{eq:spherical}), we have the parametric coordinates of points on the circle $k$ given by the angle $\varphi'$ for a fixed $\psi$:
\begin{equation}
\label{eq:circ1S2}
k(\varphi')=\begin{pmatrix}
	\sin\psi\cos\varphi'\\
	\sin\psi\sin\varphi'\\
	\cos\psi
\end{pmatrix}, \varphi'\in\langle0,2\pi).
\end{equation}
Varying the angle $\beta'$ (positions of $P$ on $c$), from equation~(\ref{eq:hopfcoordinates}) and the angle $\varphi'$ (positions of $Q$ on $k$ corresponding to distinct fibers $c$) we obtain the parametrization of a torus $\kappa$ covered by the fibers on $\mathcal{T}^3$:
\begin{equation}
\label{eq:torus1S3}
\kappa(\beta',\varphi')=\begin{pmatrix}
\cos\frac{\psi}{2}\cos(\varphi'+\beta')\\
\cos\frac{\psi}{2}\sin(\varphi'+\beta')\\
\sin\frac{\psi}{2}\cos(\beta')\\
\sin\frac{\psi}{2}\sin(\beta')
\end{pmatrix}, \beta',\varphi' \in \langle0,2\pi).
\end{equation}
Fig.~\ref{fig:hopf_torus1_circ_S} shows the double orthogonal projection of the torus $\kappa$ and its generating circles corresponding to points on the circle $k$ with their stereographic images.

The conjugated images $\kappa_1$ and $\kappa_2$ of this torus are parts of cylindrical surfaces of revolution in $\mathcal{T}_1^3$ and $\mathcal{T}_2^3$. This is a straightforward consequence of the relationship between the point $Q$ and the angle $\varphi'$. From equation~(\ref{eq:torus1S3}) with a fixed $\psi$, we obtain the $\Xi$-image $\kappa_1$ as a part of a cylindrical surface of revolution about the axis parallel to $z$ in the 3-space $\Xi(x,y,z)$. Similarly, the $\Omega$-image $\kappa_2$ is a part of a cylindrical surface of revolution about the axis parallel to $x$ in the 3-space $\Omega(x,z,w)$.

\subsection{Hopf torus of a circle with diameter parallel to $z$}
\label{ss:2}
\begin{figure*}[!htb]
\centering
\subfloat[]{\includegraphics[height=0.25\textheight]{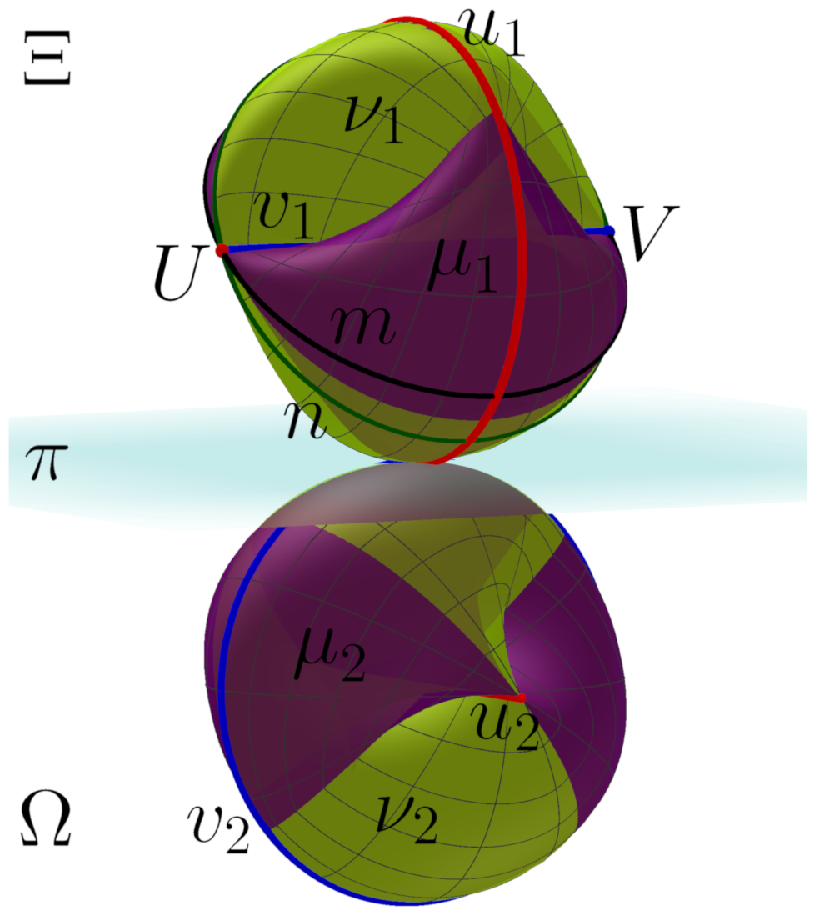}%
\label{fig:hopf_2tori_proj}}\hfil
\subfloat[]{\includegraphics[height=0.25\textheight, trim=0 20 0 20, clip]{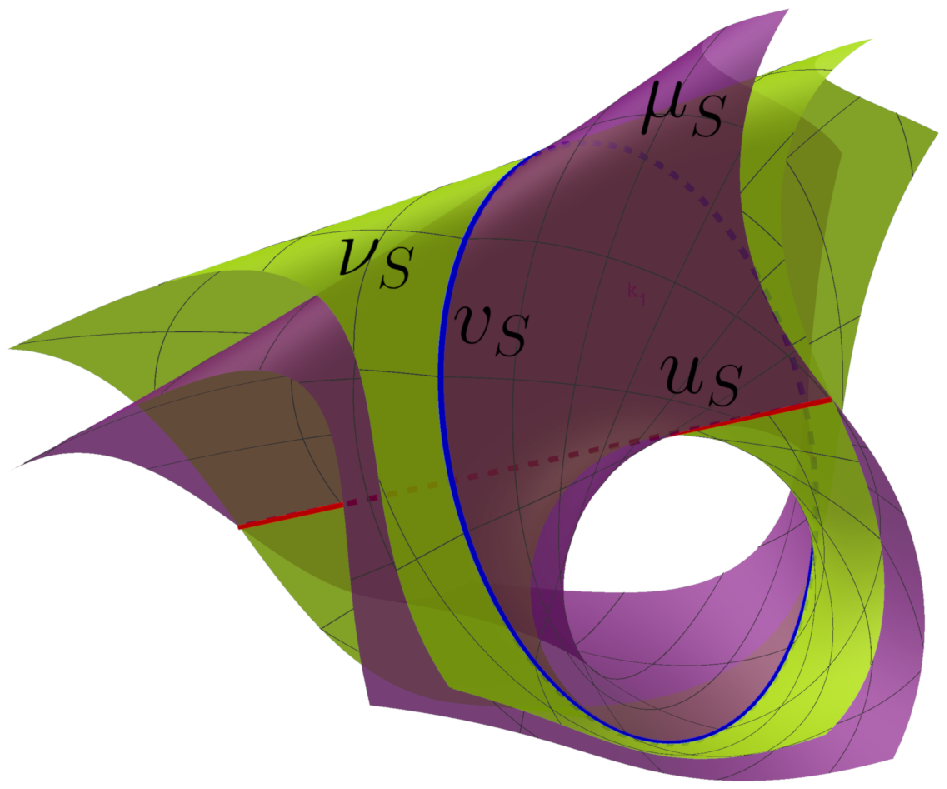}%
\label{fig:hopf_2tori_S}}
\caption{(a) Conjugated images of the tori $\mu$ and $\nu$ corresponding to the circles $m$ and $n$, respectively. The circles $m,n$ pass through the points $U,V$ on the diameter of $\mathcal{B}^2$ parallel to the $z$-axis. The fibers $u$ and $v$ above points $U$ and $V$, respectively, lie on both tori. (b) Stereographic images $\mu_S$ and $\nu_S$ of the tori. The stereographic images $u_S$ and $v_S$ of the fibers $u$ and $v$ are the intersecting line and circle of $\mu_S$ and $\nu_S$.\newline 
In the interactive model \texttt{https://www.geogebra.org/m/k94dvfpx} (or Supp. File~7), the user can manipulate the circles $m,n$ by varying the parameter $\varphi$ in equation~(\ref{eq:circ1S2+}). The projections of tori $\mu$ and $\nu$ vary dependently on $m,n$. The visibility of the objects on $\mathcal{B}^2$, the conjugated images in the double orthogonal projection, or the stereographic image can be turned off.}
\label{fig:hopf_2tori}
\end{figure*}

Let $Q$ be a point on $\mathcal{B}^2$ and $m$ a circle with a diameter parallel to $z$ through the point $Q$ (Fig.~\ref{fig:hopf_tori2}, Suppl. File~6). The circle $m$ has a parametric representation for a fixed angle $\varphi$ and variable $\psi'$ (from equation~(\ref{eq:spherical})):
\begin{equation}
\label{eq:circ2S2}
m(\psi')=\begin{pmatrix}
	\sin\psi'\cos\varphi\\
	\sin\psi'\sin\varphi\\
	\cos\psi'
\end{pmatrix}, \psi'\in\langle0,\pi\rangle.
\end{equation}
Varying the angle $\psi'$ (positions of $Q$ on $m$) changes the moduli $r'_A$ and $r'_B$ of $P$ on the corresponding fiber~$c$ (from equation~(\ref{eq:hopfcoordinates})). More precisely, $r'_A=\cos\gamma'=\cos\frac{\psi'}{2}$ and $r'_B=\sin\gamma'=\sin\frac{\psi'}{2}$ with the fixed angle $\varphi$ induce a family of non-intersecting circular fibers $c(\beta')$ generating a torus:
\begin{equation}
\label{eq:torus2S3}
\mu(\beta',\psi')=\begin{pmatrix}
\cos\frac{\psi'}{2}\cos(\varphi+\beta')\\
\cos\frac{\psi'}{2}\sin(\varphi+\beta')\\
\sin\frac{\psi'}{2}\cos(\beta')\\
\sin\frac{\psi'}{2}\sin(\beta')
\end{pmatrix}, \beta'\in \langle0,2\pi),\psi'\in\langle0,\pi\rangle.
\end{equation}
In contrast to the previous case, in which the torus $\kappa$ was generated by the parameter $\varphi'$, on which the first two coordinates depend, the torus $\mu$ is generated by the parameter $\psi'$, on which all four coordinates depend. Therefore, the conjugated images of the torus $\mu$ are more twisted. For the choice $\psi'=0$ and $\beta'=\frac{\pi}2-\varphi$, we always obtain the point $[0,1,0,0]$ lying on the torus $\mu$. This is the center $N$ of the stereographic projection before applying the translation for visualization purposes, so the torus $\mu$ always contains a point that is stereographically projected to infinity. Moreover, it contains a fiber through this point, too. The stereographic image of this fiber is a line. See Fig.~\ref{fig:hopf_torus2_circ_S} for the full illustration of the double orthogonal projection of the torus covered by its circles and their stereographic images. Let $\mu$ and $\nu$ be tori generated by great circles $m$ and $n$ with a diameter parallel to the $z$-axis (Fig.~\ref{fig:hopf_2tori}). The circles $m$ and $n$ intersect in antipodal points $U$ and $W$, which correspond to fibers $u$ and $v$, respectively. These fibers were depicted earlier in Figs.~$\ref{fig:hopf_fiber_positions}$ (a) and~(b). Therefore, the tori $\mu$ and $\nu$ have two common fibers such that one conjugated image is always a segment and the second image is a great circle. By the above-mentioned property of the fiber through the center of stereographic projection ($u$ in Fig.~\ref{fig:hopf_2tori}), the stereographic images $\mu_S$ and $\nu_S$ of the tori have a common line $u_S$.

\subsection{Nested Hopf tori corresponding to families of circles on $\mathcal{B}^2$}
\begin{figure*}[!]
\centering
\subfloat[]{\includegraphics[width=\linewidth, trim=0 0 0 0, clip]{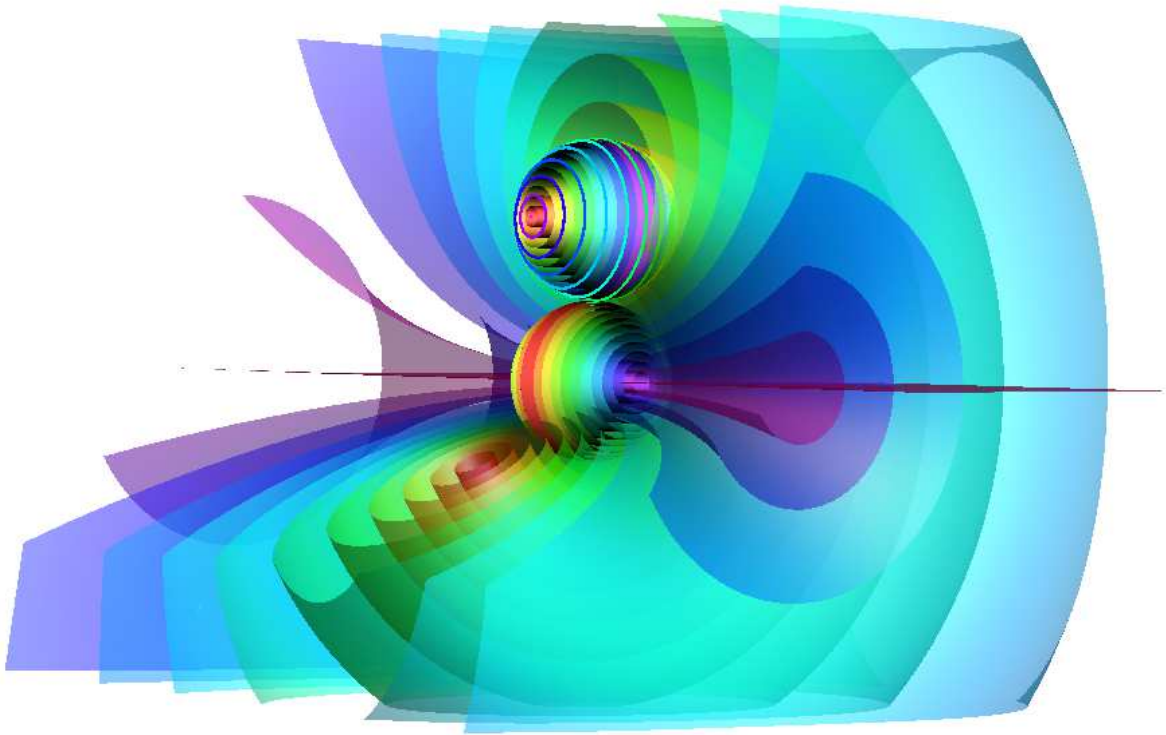} %
\label{fig:nested_tori1_S}}\\\bigskip
\subfloat[]{\includegraphics[width=0.9\linewidth]{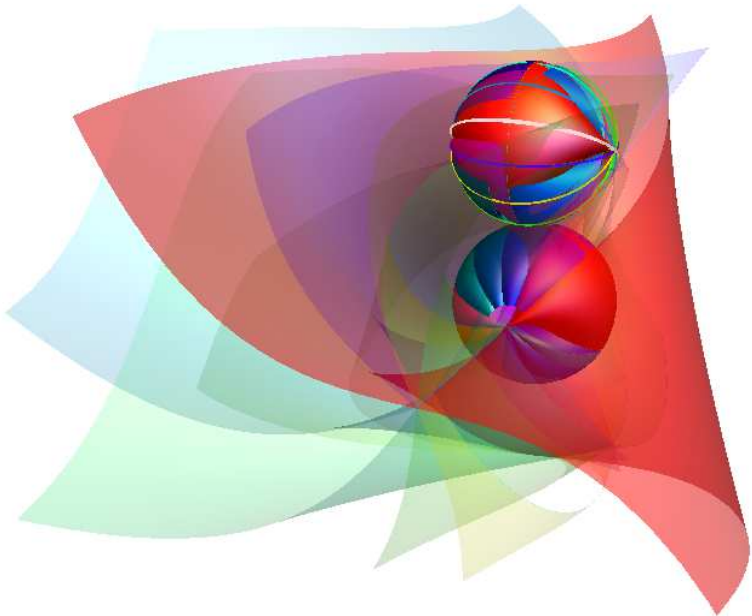}%
\label{fig:nested_tori2_S}}
\caption{(a) A family of circles on $\mathcal{B}^2$ parallel to $(x,y)$  and with a diameter parallel to $z$; (b) the corresponding nested tori on $\mathcal{T}^3$, and their stereographic projection onto $\Omega(x,z,w)$. Colors (shades) refer to mutually related objects; in (b) the images of the torus highlighted in red corresponds to the white great circle on $\mathcal{B}^2$. The visualizations are based on the parametrization in equation~(\ref{eq:omega+}).}
\label{fig:nested_tori_S}
\end{figure*}
We summarize the toroidal structure of a 3-sphere in the following visualizations. For each circle $k$ parallel to the $(x,y)$-plane on the \mbox{2-sphere} $\mathcal{B}^2$ we obtain a torus $\kappa$. Fig.~\ref{fig:nested_tori1_S} gives a model (see Suppl. File~8 for a video animation) of nested tori $\kappa$ on $\mathcal{T}^3$ corresponding to circles $k$ on $\mathcal{B}^2$. This family of disjoint tori contains only one fiber through the center of the stereographic projection, and hence the tori appear in the stereographic projection as nested tori of revolution including one line, which is their axis.

In the second case (Fig.~\ref{fig:nested_tori2_S}, Suppl. File~9), the family of circles $m$ on $\mathcal{B}^2$ with a diameter parallel to the $z$-axis forms nested tori $\mu$ on $\mathcal{T}^3$. Each of these tori contains two common fibers in the special positions shown in Fig.~\ref{fig:hopf_2tori}. The stereographic image of one of the fibers is a line, which is the common line for all the stereographic images of the tori.

 These families of tori cover the \mbox{3-sphere} $\mathcal{T}^3$ reparametrized by variables $\beta',\varphi'$ and~$\psi'$ as:
\begin{equation}
\begin{split}
\label{eq:omega}
\mathcal{T}^3(\beta',\varphi',\psi')=\begin{pmatrix}
\cos\frac{\psi'}{2}\cos(\varphi'+\beta')\\
\cos\frac{\psi'}{2}\sin(\varphi'+\beta')\\
\sin\frac{\psi'}{2}\cos(\beta')\\
\sin\frac{\psi'}{2}\sin(\beta')
\end{pmatrix}, \\ \beta',\varphi' \in \langle0,2\pi), \psi'\in \langle0,\pi\rangle.
\end{split}
\end{equation}

{
\section{Further applications}
\label{sec:applications}
In the last section, we provide a few brief references to geometrically challenging applications of the Hopf fibration.
\subsection{Cyclic surfaces}
\begin{figure*}[!htb]
\centering
\subfloat[]{\includegraphics[height=0.25\textheight]{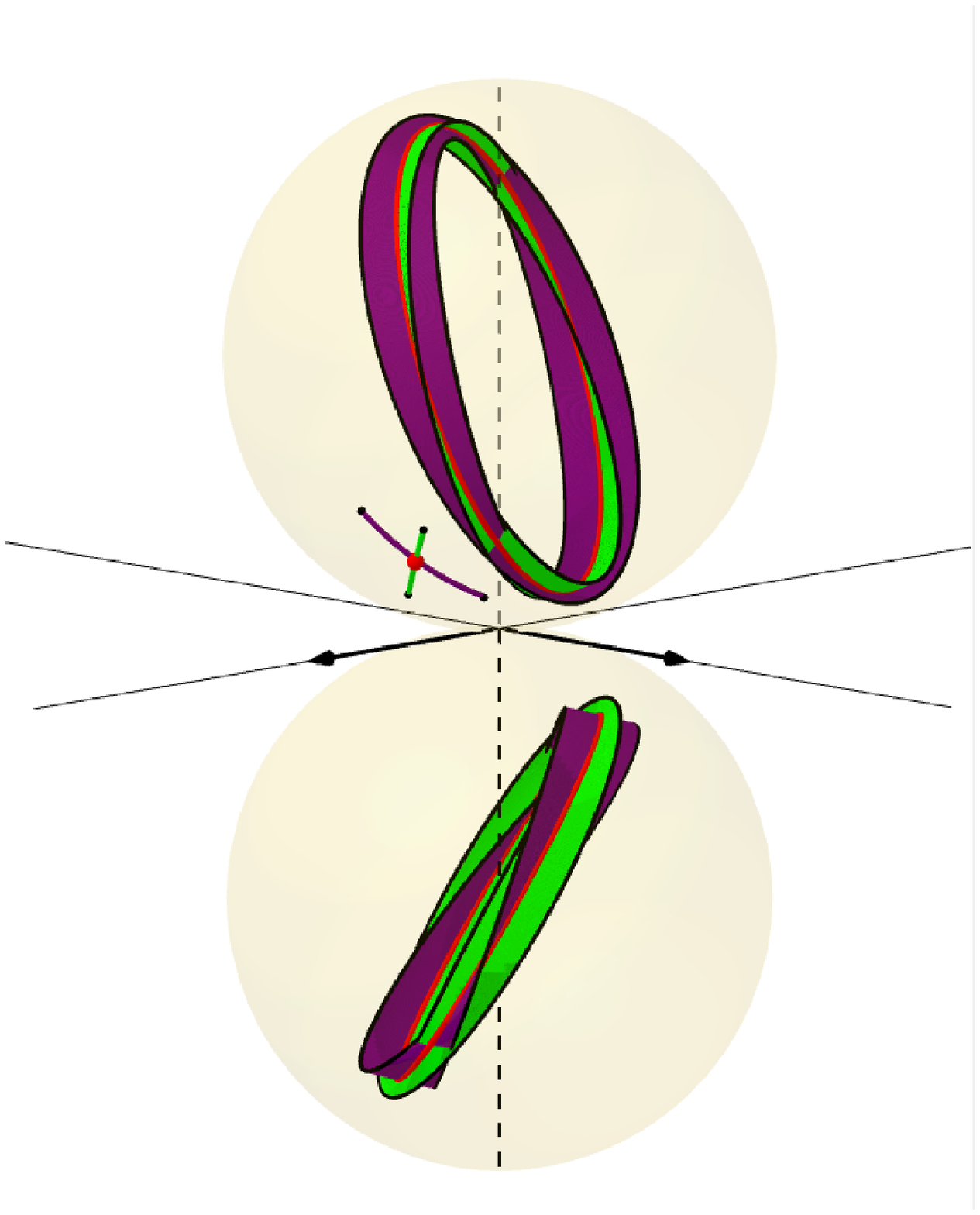}%
\label{fig:hopf_intersection_proj}}\hfil
\subfloat[]{\includegraphics[height=0.25\textheight]{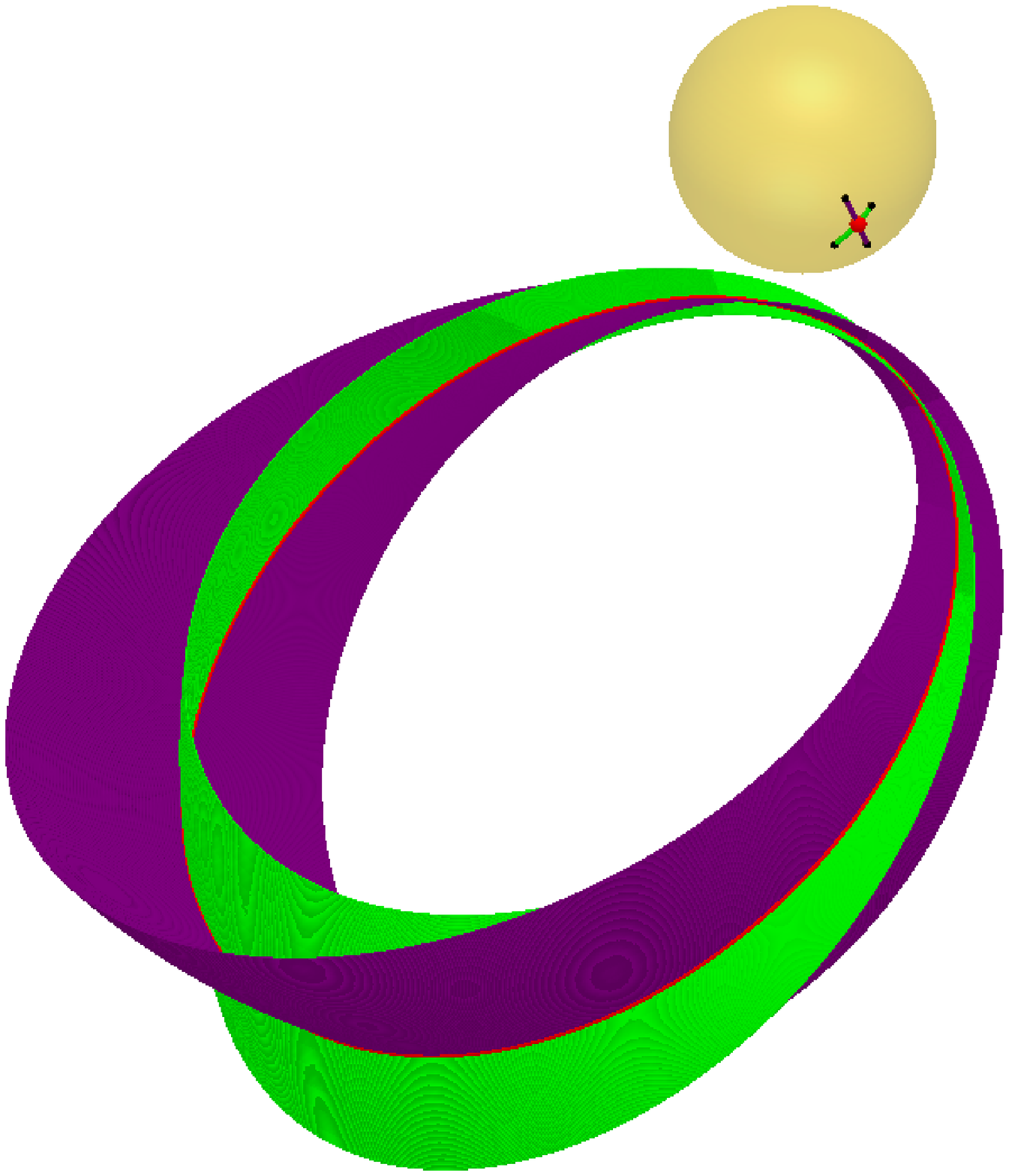}%
\label{fig:hopf_intersection_S}}
\caption{The intersection point of two curves on a 2-sphere $\mathcal{B}^2$. The corresponding fibers in the Hopf fibration form two cyclic surfaces on the 3-sphere $\mathcal{T}^3$ intersecting in one circle. The situation is visualized in the double orthogonal projection (a) and stereographic projection (b).}
\label{fig:hopf_intersection}
\end{figure*}

Using the above-mentioned method of visualization and construction of the Hopf fibration, we can study related properties of geometric surfaces and design shapes that are formed by disjoint circles. If we consider a point moving along an arbitrary curve on the 2-sphere $\mathcal{B}^2$, the motion of the corresponding Hopf fiber creates a cyclic surface consisting of disjoint circles of variable radius on the 3-sphere $\mathcal{T}^3$ embedded in $\mathbb{R}^4$. Consequently, two curves intersecting in one point on a 2-sphere create two cyclic surfaces with only one common circle (the Hopf fiber of the point of intersection) (Fig.~\ref{fig:hopf_intersection_proj}). Stereographic projection preserves circles (up to a circle through the center of projection) and so the stereographic images are cyclic surfaces intersecting in a circle in $\mathbb{R}^3$, too (Fig.~\ref{fig:hopf_intersection_S}).

\begin{figure*}[!htb]
\centering
\subfloat[]{\includegraphics[height=0.25\textheight, trim=100 0 100 0, clip]{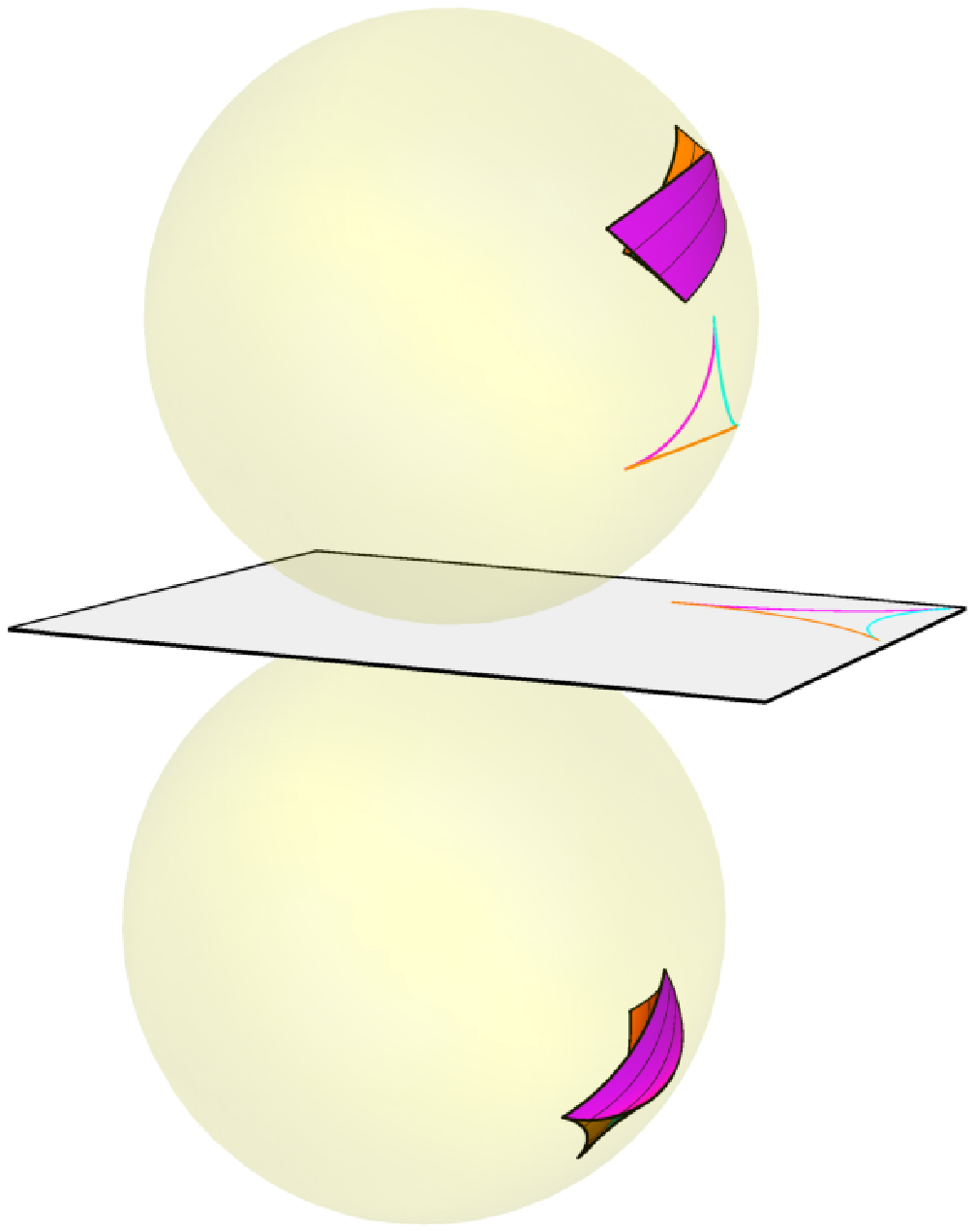}%
\label{fig:hopf_arcs_proj}}\hfil
\subfloat[]{\includegraphics[height=0.25\textheight, trim=100 100 100 100, clip]{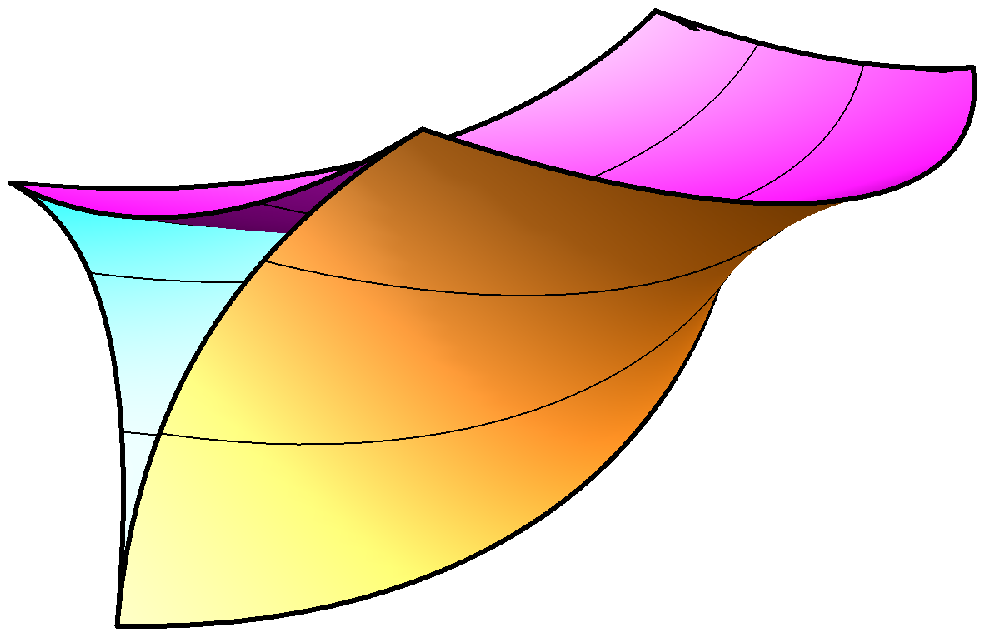}%
\label{fig:hopf_arcs_S}}\\
\subfloat[]{\includegraphics[height=0.25\textheight, trim=100 0 100 0, clip]{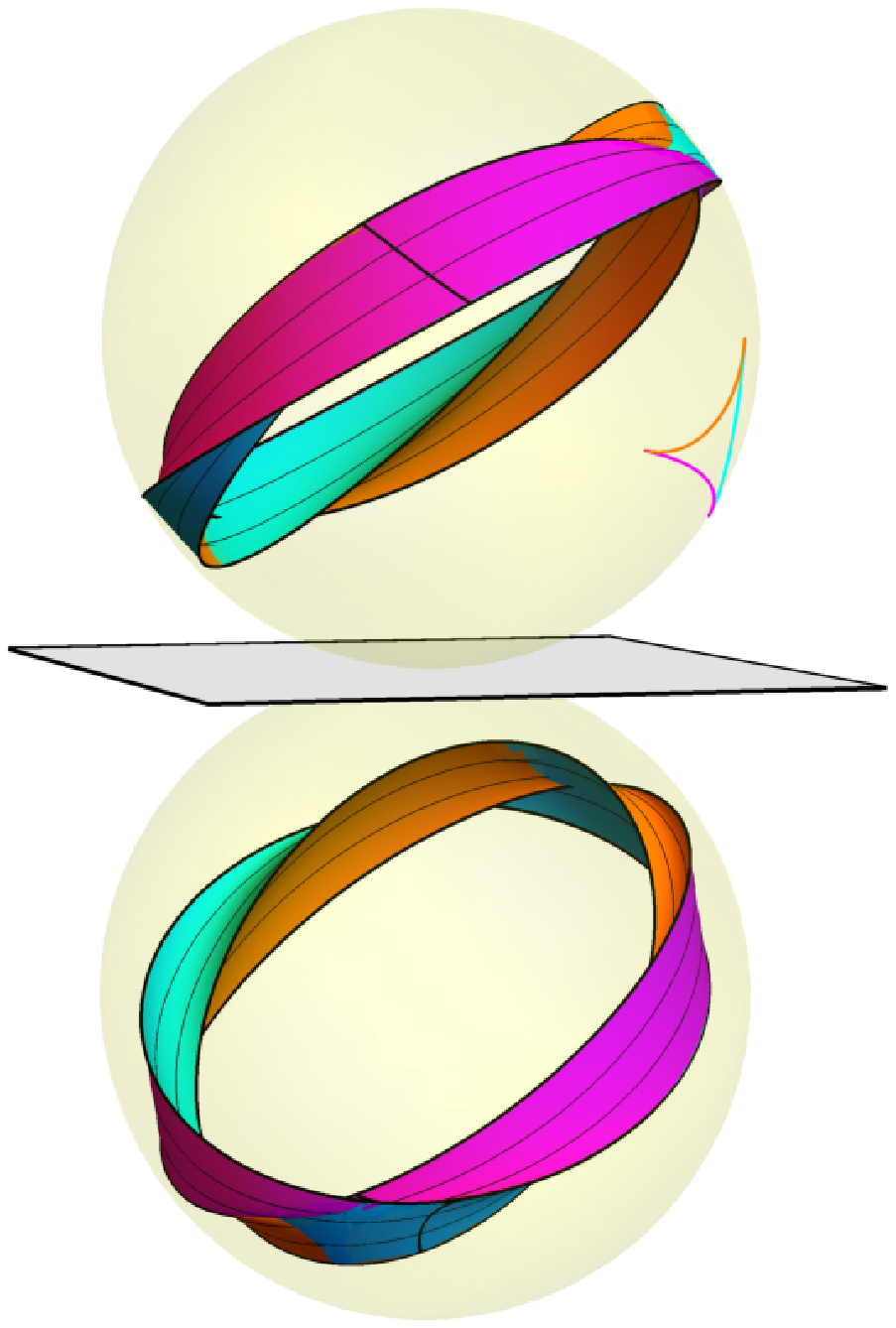}%
\label{fig:hopf_arcsfull_proj}}\hfil
\subfloat[]{\includegraphics[height=0.25\textheight, trim=100 100 100 100, clip]{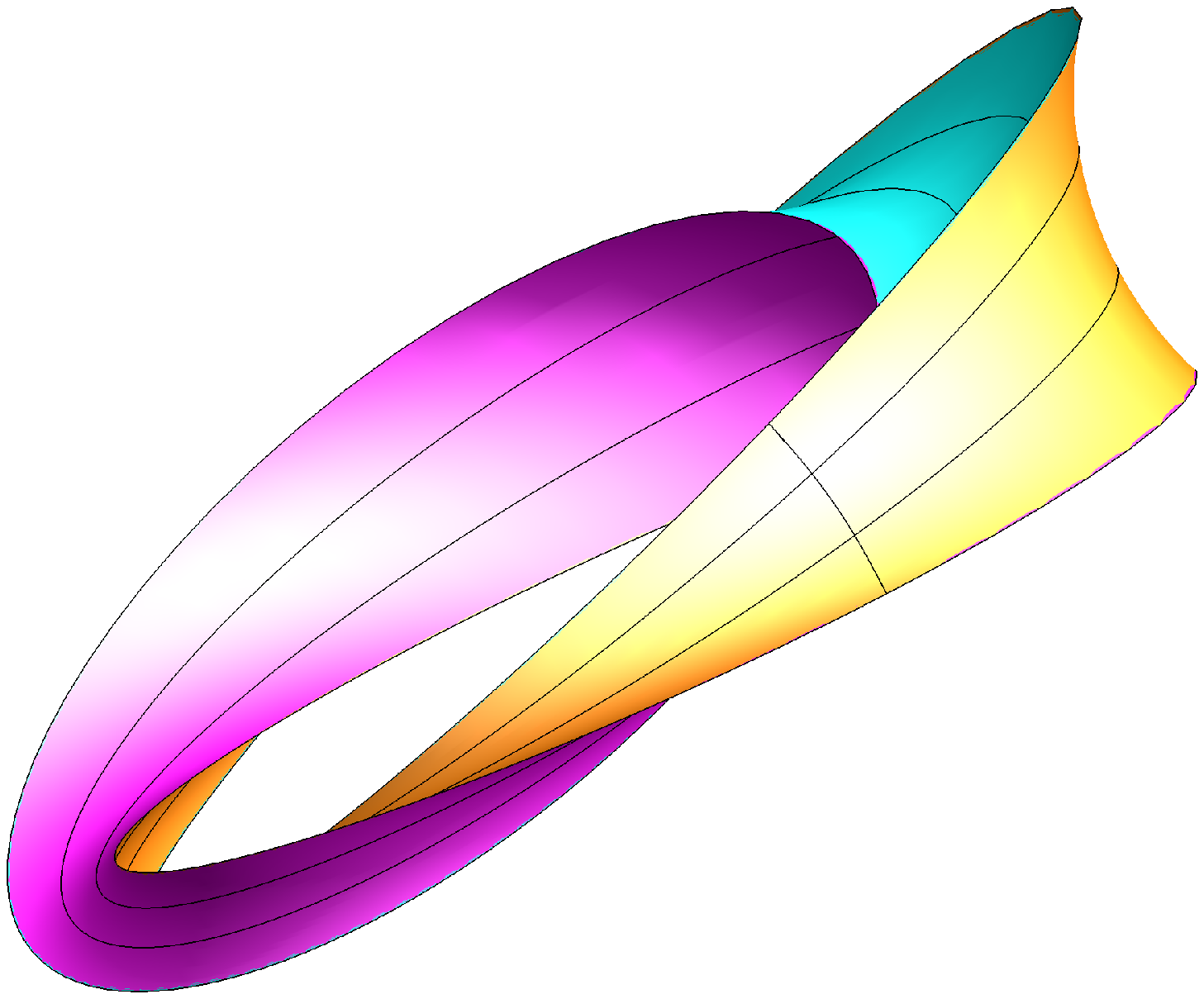}%
\label{fig:hopf_arcsfull_S}}
\caption{A planar shape created by three circular arcs in the $(x,z)$-plane is stereographically projected onto the 2-sphere $\mathcal{B}^2\in\mathbb{R}^3$. The corresponding fibers in the Hopf fibration into $\mathcal{T}^3\in\mathbb{R}^4$ form three connected parts of cyclic surfaces, each two connected by a common circle. A part of the shape is depicted in (a) the double orthogonal projection and (b) the stereographic projection. The whole shape is visualized in (c) the double orthogonal projection and (d) the stereographic projection.}
\label{fig:hopf_arcs}
\end{figure*}

Finally, we can construct orthogonal and stereographic images of shapes consisting of cyclic surfaces or their parts connected by common circles. The case in Fig.~\ref{fig:hopf_arcs} shows a union of three circular arcs (the vertices are tangent points of the corresponding circles) stereographically projected from the $(x,z)$-plane onto a 2-sphere $\mathcal{B}^2$. Then, we apply the inverse Hopf projection and construct conjugated images of the corresponding surfaces in $\mathcal{T}^3$ in the double orthogonal projection. After the stereographic projection from $\mathcal{T}^3$ to the 3-space $\Omega(x,z,w)$ we obtain a three-dimensional model of the shape as the union of parts of cyclic surfaces. The common points of each pair of circular arcs become the common circles of each pair of parts of the cyclic surfaces.

\subsection{Four-dimensional modulations}
The Hopf fibration was used as a constructive tool in classical optical communications to design four-dimensional modulations in \cite{Rodrigues2018}. From the geometric point of view, $n$PolSK-$m$PSK modulations are constructed by $n$ vertices of a polyhedron inscribed in the base 2-sphere generating $n$ fibers in the 3-sphere, and each fiber contains $m$ points. The authors demonstrated their main results on 14PolSK-8PSK modulation generated by tetrakis hexahedron, which is a union of a hexahedron and octahedron with 14 vertices. The visualization of this arrangement supplemented with its double orthogonal projection is in Fig.~\ref{fig:modulation}.
\begin{figure}[!htb]
\centering
\includegraphics[width=.85\textwidth]{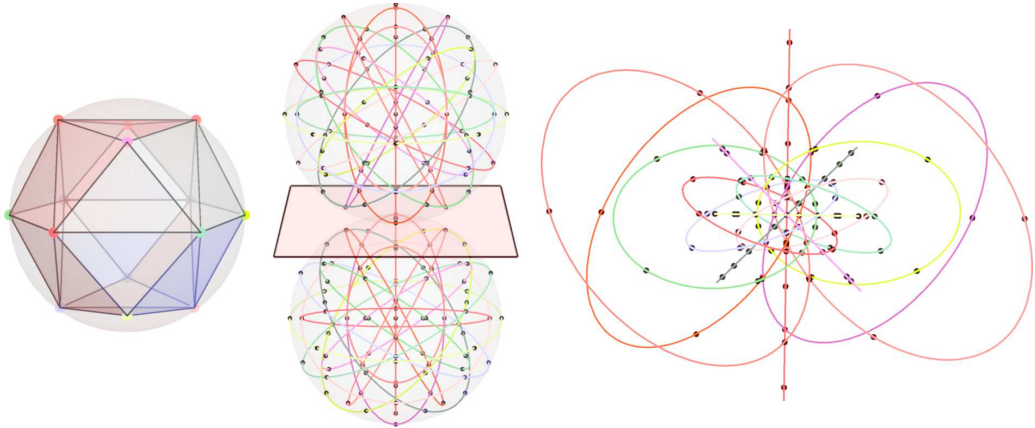}
\caption{Visualization of the 14PolSK-8PSK modulation in the double orthogonal projection and stereographic projection given by a tetrakis hexahedron with 14 vertices.}
\label{fig:modulation}
\end{figure}

\subsection{Twisted filaments}
\begin{figure}[!htb]
	\centering
	\includegraphics[width=.85\linewidth]{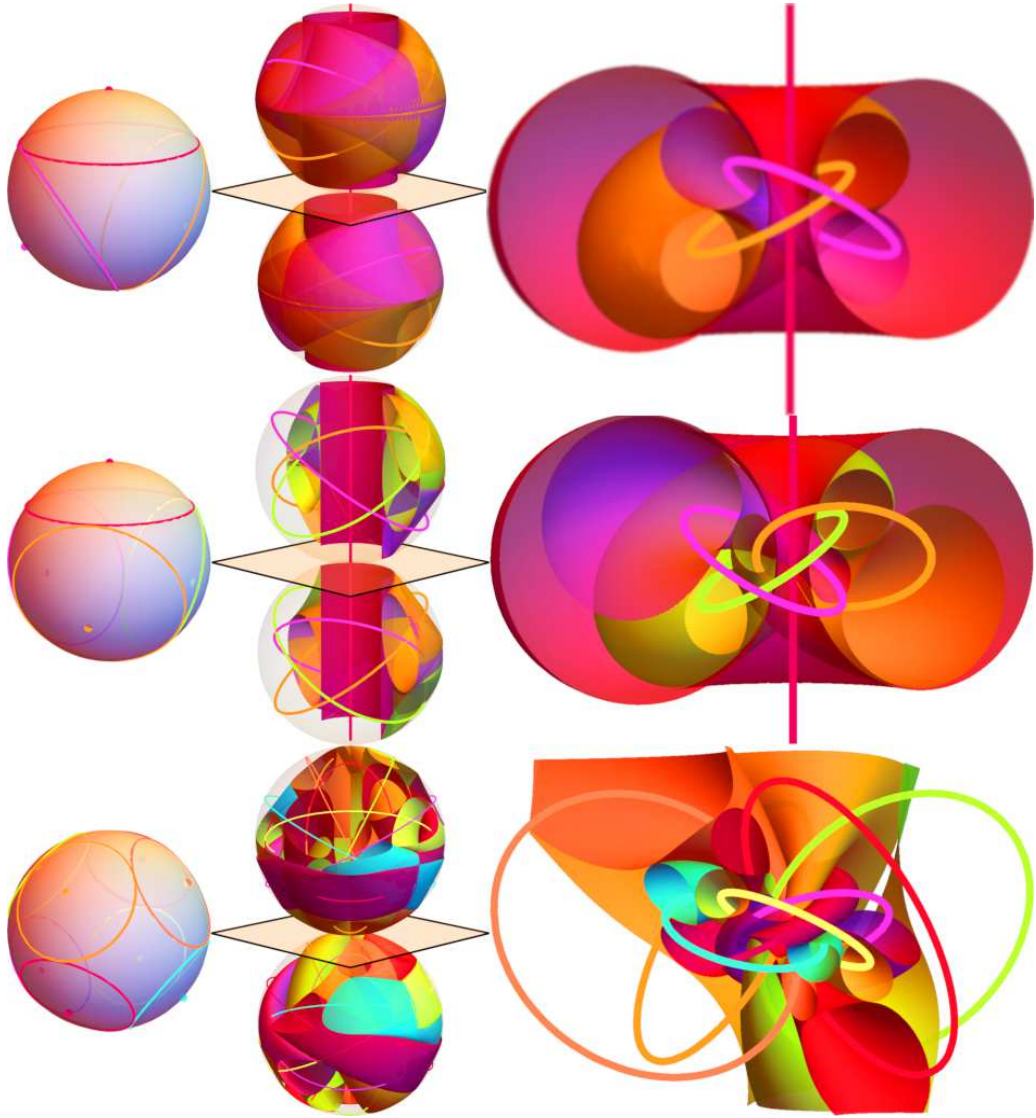}%
\caption{Triangular, tetrahedral, and hexahedral arrangements of vertices on the base 2-sphere and their corresponding twisted filaments with backbones visualized in the double orthogonal projection and in the stereographic projection.}
\label{fig:packings1}
\end{figure}

Our final application is inspired by twisted toroidal structures appearing in biological and synthetic materials. In this context, \cite{Atkinson2019, Grason2015, Kleman1985} studied the problem of twisted filament packings. Geometric models are constructed with the use of the Hopf fibration, in which fibers play the role of filament backbones. To construct equally spaced filaments in the 3-sphere $\mathcal{T}^3$, we consider equally spaced disks on the base 2-sphere $\mathcal{B}^2$. Visualizations of such arrangements based on the vertices of Platonic solids (and triangular case) projected to the 2-sphere $\mathcal{B}^2$ in the double orthogonal projection and also stereographic projection are in Figs.~\ref{fig:packings1} and \ref{fig:packings2}. Tangent points of the circles on $\mathcal{B}^2$ correspond to tangent circles along the filaments in $\mathcal{T}^3$. Consequently, the number of neighboring filaments in $\mathcal{T}^3$ is the same as the number of neighboring disks in~$\mathcal{B}^2$. 

\begin{figure}[!htb]
	\centering
	\includegraphics[width=.85\linewidth]{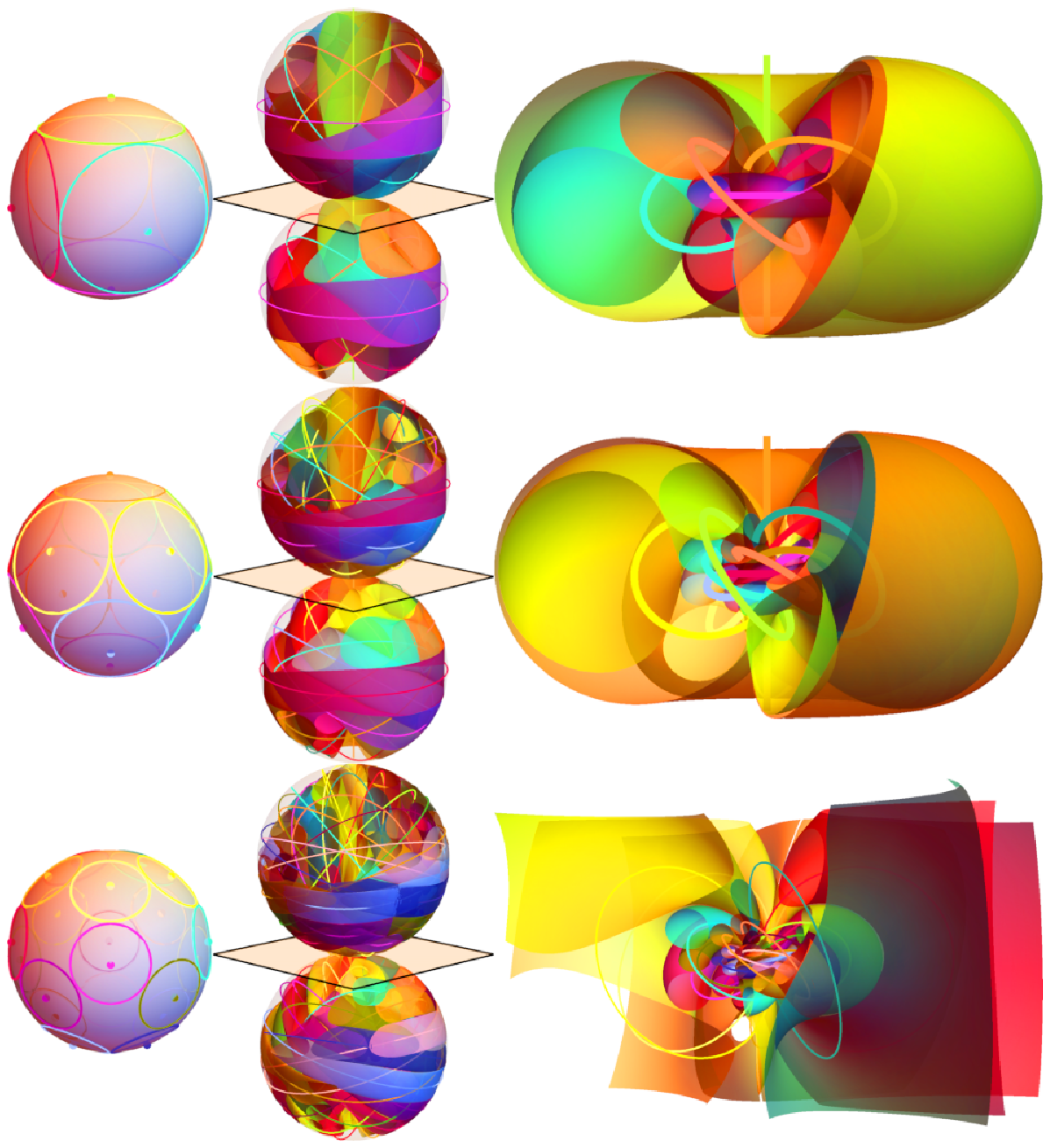}%
\caption{Octahedral, icosahedral, and dodecahedral arrangements of vertices on the base 2-sphere and their corresponding twisted filaments with backbones visualized in the double orthogonal projection and in the stereographic projection.}
\label{fig:packings2}
\end{figure}

As an example of more complex structure of filaments packing, we chose buckminsterfullerene with 60 vertices, 12 pentagonal and 20 hexagonal faces (Fig.~\ref{fig:packingsful}).
\begin{figure}[!htb]
\centering
\includegraphics[width=.65\textwidth]{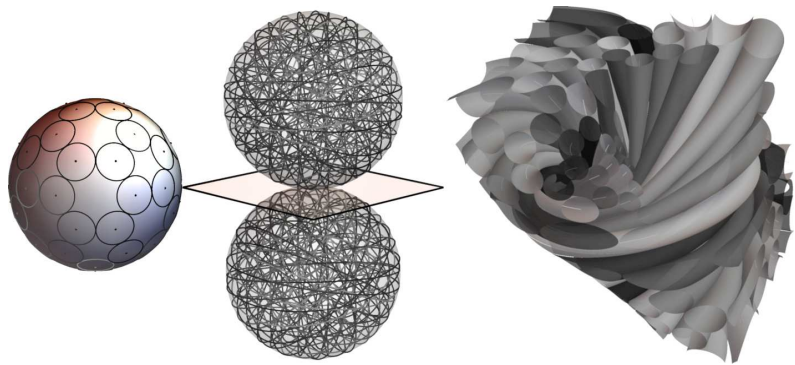}
\caption{Vertices of a buckminsterfullerene projected to the base 2-sphere, their corresponding filament backbones in the double-orthogonal projection, and a close-up to stereographic images of the corresponding twisted filaments.}
\label{fig:packingsful}
\end{figure}

\section{Conclusion}
By the use of elementary constructive tools in the double orthogonal projection of \mbox{4-space} onto two mutually perpendicular \mbox{3-spaces} rotated into one 3-dimensional modeling space, we have described a synthetic step-by-step construction of a Hopf fiber on a \mbox{3-sphere} embedded in \mbox{4-space} that corresponds to a point on a \mbox{2-sphere}. The virtual modeling space is accessible in supplementary interactive models created in the interactive 3D geometric software {\sl GeoGebra 5}, in which the reader can intuitively manipulate fundamental objects and achieve a sense of the fourth dimension through two interlinked 3-dimensional models. The choice of the method of visualization plays a significant role in several aspects. First, two conjugated three-dimensional images of a four-dimensional object carry all the necessary information to determine this object uniquely. For example, if we had only one image of parallel sections of a 3-sphere in Fig.~\ref{fig:2sphere_sections}, we would miss important details for reconstruction. Interpretation of two three-dimensional images as one object indeed assumes some training and experience. However, in the case of projection onto a plane, we would need at least three images for a visual representation. Another advantage of the double orthogonal projection is that synthetic constructions generalize constructions in Monge's projection. Thus, for example, the localization of a point on a 3-sphere or constructions of stereographic images are elementary.
Furthermore, with the use of this technique, we have visualized a torus formed by the Hopf fibers corresponding to a circle on the \mbox{2-sphere} and shown it in two different positions of the circles on the \mbox{2-sphere}. The Hopf fibration is usually visualized in the stereographic projection based on an analytic representation. Since intuition in four-dimensional visualization, including the double orthogonal projection, is often misleading, stereographic images, in this case, support our reasoning. Therefore, the tori were projected to synthetically constructed stereographic images in the modeling \mbox{3-space}, revealing the true nature of the Hopf fibration. The final visualizations show how the points of the \mbox{2-sphere} cover the \mbox{3-sphere} by their corresponding fibers on the nested tori, and the whole modeling \mbox{3-space} when stereographically projected. In this way we have built a mathematical visualization of the Hopf fibration in which we have presented a constructive connection between the base space (a 2-sphere), fibers covering the total space (a 3-sphere), and the stereographic images, and through which we can study and explain properties of the Hopf fibration in a way that does not depend on its analytic description. However, for the purposes of verification and implementation, the objects are supported by their parametric representations used in the classical analytic approach.

We finished by giving a short application of the Hopf fibration for constructing shapes on the 3-sphere and in stereographic projection in 3-space. The method provided promises future applications in visualizing other curves on a \mbox{2-sphere} and their corresponding surfaces generated by their Hopf fibers on a \mbox{3-sphere} in \mbox{4-space}, along with their corresponding stereographic images in the modeling \mbox{3-space}. The double orthogonal projection method is also likely to be useful for visualizing and analyzing the properties of further 3-manifolds embedded in \mbox{4-space}.

\section{Acknowledgment}
I wish to thank Lada Peksová (Charles University) for her help with the topological aspects and Filip Beran (Charles University) for revision and valuable suggestions. I also thank the anonymous referees for their valuable suggestions, which helped to improve the exposition. The work was supported by the grant Horizon 2020 MaTeK Enhancement of research excellence in Mathematics Teacher Knowledge (No 951822).
\bibliographystyle{apa}
\bibliography{zamboj_manuscript_jcde_rev3_clean}

\newpage
\appendix{
\section{Parametrizations relevant to the figures}
{\label{appendix}
Parametrizations in $\mathbb{R}^4$ are given with $(x,y,z,w)$-coordinates. For the implementation of visualizations in the double orthogonal projection in the modeling space with $(x,y,z)$-coordinates, we decompose the images so that a $\Xi$-image has $(x,y,z)$-coordinates and an \mbox{$\Omega$-image} has $(x,-w,z)$-coordinates.
\subsubsection*{Figs.~\ref{fig:hopf_SBS1} and \ref{fig:hopf_SBS}}
A point $Q$ on the \mbox{2-sphere} $\mathcal{B}^2$ (cf. Equation~\ref{eq:spherical}):
\begin{equation}
\begin{split}
\label{eq:parS2+}
Q  &=\begin{pmatrix}
	\sin\psi\cos\varphi\\
	\sin\psi\sin\varphi +1\\
	\cos\psi
\end{pmatrix},\\ \psi &\in\langle0,\pi\rangle,\varphi\in\langle0,2\pi).
\end{split}
\end{equation} 
Further on $\psi$ and $\varphi$ are fixed, and $\gamma=\frac{\psi}{2},\varphi=\alpha-\beta$, for $\beta \in \langle0,2\pi)$.\\
A point $P$ on the \mbox{3-sphere} $\mathcal{T}^3$ (cf. Equation~\ref{eq:hopfcoordinates}):
\begin{equation}
\label{eq:parS3gon+}
P=\begin{pmatrix}
\cos\frac{\psi}{2}\cos(\varphi+\beta)\\
\cos\frac{\psi}{2}\sin(\varphi+\beta)+1\\
\sin\frac{\psi}{2}\cos\beta\\
\sin\frac{\psi}{2}\sin\beta+1
\end{pmatrix},
\end{equation}
and also the parametrization of the circle $c(\beta)$ by the variable~$\beta$.\\
The point $P_S$ in $\Omega(x,z,w)$ -- the stereographic image of the point~$P$ from the center $N=[0,2,0,1]$
\begin{equation}
\label{eq:parS3gonSter+}
P_S=\begin{pmatrix}
\displaystyle\frac{2\cos\frac{\psi}{2} \cos(\varphi + \beta)}{1 - \cos\frac{\psi}{2}\sin(\varphi+ \beta)}\\
\displaystyle\frac{2\sin\frac{\psi}{2}\cos\beta}{1 - \cos\frac{\psi}{2}\sin(\varphi+\beta)}\\
\displaystyle\frac{2\sin\frac{\psi}{2}\sin\beta}{1 - \cos\frac{\psi}{2}\sin(\varphi+\beta)} + 1
\end{pmatrix},
\end{equation}
and also the parametrization of the circle $c_S(\beta)$ by the variable~$\beta$.
\subsubsection*{Figs.~\ref{fig:hopf_tori} and \ref{fig:hopf_tori_circ_S}}
A circle $k$ on $\mathcal{B}^2$ for a fixed $\psi$ (cf. Equation~\ref{eq:circ1S2}):
\begin{equation}
\label{eq:circ1S2+}
k(\varphi')=\begin{pmatrix}
	\sin\psi\cos\varphi'\\
	\sin\psi\sin\varphi'\\
	\cos\psi+1
\end{pmatrix}, \varphi'\in\langle0,2\pi).
\end{equation}
The torus $\kappa$ corresponding to the circle $k$ (cf. Equation~\ref{eq:torus1S3}):
\begin{equation}
\label{eq:torus1S3+}
\kappa(\beta',\varphi')=\begin{pmatrix}
\cos\frac{\psi}{2}\cos(\varphi'+\beta')\\
\cos\frac{\psi}{2}\sin(\varphi'+\beta')+1\\
\sin\frac{\psi}{2}\cos(\beta')\\
\sin\frac{\psi}{2}\sin(\beta')+1
\end{pmatrix}, \beta',\varphi' \in \langle0,2\pi).
\end{equation}
The stereographic image of the torus $\kappa$ in $\Omega(x,z,w)$ (cf. Equation~\ref{eq:omega}):
\begin{equation}
\label{eq:torus1ster+}
\kappa(\beta',\varphi')=\begin{pmatrix}
\displaystyle\frac{2\cos\frac{\psi}{2} \cos(\varphi' + \beta') }{1 - \cos\frac{\psi}{2} \sin(\varphi' + \beta')}\\
\displaystyle\frac{2\sin\frac{\psi}{2} \cos(\beta')}{1 - \cos\frac{\psi}{2} \sin(\varphi' + \beta')}\\
\displaystyle\frac{2\sin\frac{\psi}{2} \sin(\beta')}{1 - \cos\frac{\psi}{2} \sin(\varphi' + \beta')} + 1
\end{pmatrix},\beta',\varphi'\in \langle0,2\pi).
\end{equation}
Analogously, a circle $m$ and torus $\mu$, along with their stereographic images, differ from $k$ and $\kappa$ only by fixing $\varphi$ and making the variable $\psi'$. 
\subsubsection*{Fig.~\ref{fig:nested_tori_S}}
The stereographic images of the nested tori in $\Omega(x,z,w)$ have the following parametric representation (cf. Equation~\ref{eq:omega}):
\begin{equation}
\begin{split}
\label{eq:omega+}
\mathcal{T}^3(\beta',\varphi',\psi')&=\begin{pmatrix}
\displaystyle\frac{2\cos\frac{\psi'}{2} \cos(\varphi' + \beta') }{1 - \cos\frac{\psi'}{2} \sin(\varphi' + \beta')}\vspace{1pt}\\
\displaystyle\frac{2\sin\frac{\psi'}{2} \cos(\beta')}{1 - \cos\frac{\psi'}{2} \sin(\varphi' + \beta')}\vspace{1pt}\\
\displaystyle\frac{2\sin\frac{\psi'}{2} \sin(\beta')}{1 - \cos\frac{\psi'}{2} \sin(\varphi' + \beta')} + 1
\end{pmatrix}.
\end{split}
\end{equation}
The stereographic image in Fig.~\ref{fig:nested_tori1_S} is for the sake of clarity restricted to $\beta'\in\langle\frac{\pi}{6}, \frac{3\pi}{2}\rangle$, and $\psi'=k\frac{\pi}{12}$ for $k\in\{0,1,\dots,12\}$ is chosen to be the leading variable. In Fig.~\ref{fig:nested_tori2_S} the $\Xi$ and $\Omega$-images $\mu_1$ and $\mu_2$ of the tori $\mu$ are, apart of the one highlighted in red, restricted to $\beta'\in\langle\frac{\pi}{3},\frac{5\pi}{3}\rangle, \psi'\in\langle\frac{\pi}{6}, 2\pi\rangle$, the stereographic images are restricted to $\beta'\in\langle\frac{\pi}{6}, 11\frac{\pi}{6}\rangle$, and the leading variable is $\varphi'=k\frac{\pi}{6}$ for $k\in\{0,1,\dots,6\}$.
}

\end{document}